\documentclass{amsart}

\usepackage{url}
\usepackage{algorithm}
\usepackage{algpseudocode}
\usepackage{mathrsfs}  
\usepackage{subcaption}
\usepackage{graphicx}
\usepackage[giveninits=true]{biblatex} 
\newcommand{\norm}[1]{\left\|#1\right\|}
\newcommand{\normdg}[1]{{\left\vert\kern-0.25ex\left\vert\kern-0.25ex\left\vert #1 
    \right\vert\kern-0.25ex\right\vert\kern-0.25ex\right\vert}}

\newcommand{\jmp}[1]{[\![#1]\!]}
\newcommand{\avg}[1]{\{\!\!\{#1\}\!\!\}}

\DeclareMathOperator*{\argmin}{argmin}
\newcommand{\btensor}[1]{\underline{\boldsymbol{#1}}}
\newcommand{\bvector}[1]{\boldsymbol{#1}}

\algrenewcommand\algorithmiccomment[2][\footnotesize]{{#1\hfill\(\triangleright\) \textit{#2}}}

\newcommand{\DG}{\mathsf{DG}}
\newcommand{\VEM}{\mathsf{VEM}}
\newcommand{\D}{\mathsf{D}}
\newcommand{\E}{\mathsf{E}}
\newcommand{\Pgrad}[1]{\Pi^{\,\nabla}_{#1}}

\usepackage{xcolor}

\newtheorem{remark}[section]{Remark}
\newtheorem{definition}[section]{Definition}

\begin{document}
\title[Deep Learning Accelerated 
AMG Methods for Polytopal Discretizations]{Deep Learning Accelerated 
Algebraic Multigrid Methods for Polytopal Discretizations of Second-Order Differential Problems} 
\author[P.F. Antonietti]{Paola F. Antonietti}
\author[M. Caldana]{Matteo Caldana}
\author[L. Gentile]{Lorenzo Gentile}
\author[M. Verani]{Marco Verani}
\date{\today}
\address{MOX, Dipartimento di Matematica, Politecnico di Milano, Piazza Leonardo da Vinci 32, 20133 Milano, Italy.}
\email{paola.antonietti@polimi.it}
\email{matteo.caldana@polimi.it}
\email{lorenzo3.gentile@mail.polimi.it}
\email{marco.verani@polimi.it}

\begin{abstract}
Algebraic Multigrid (AMG) methods are state-of-the-art algebraic solvers for Partial Differential Equations. Still, their efficiency depends heavily on the choice of suitable parameters and/or ingredients. Paradigmatic examples include the so-called strong threshold parameter, which controls the algebraic coarse-grid hierarchy, as well as the smoother, i.e., the relaxation methods used on the fine grid to damp out high-frequency components of the error. In AMG, since the coarse grids are constructed algebraically (without geometric intuition), the smoother's performance is even more critical.
For the linear systems stemming from polytopal discretizations, such as Polytopal Discontinuous Galerkin (PolyDG) and Virtual Element (VEM) Methods, AMG sensitivity to such choices is even more critical due to the significant variability of the underlying meshes, which results in algebraic systems with different sparsity patterns.
In this paper, we focus on the linear systems of equations stemming from polytopal discretizations of second-order elliptic problems. We propose a novel deep learning approach that automatically tunes the strong threshold parameter and the smoother choice in AMG solvers, thereby maximizing AMG performance. 
We test various differential problems in both two- and three-dimensional settings, with heterogeneous coefficients and polygonal/polyhedral meshes, and demonstrate that the proposed approach generalizes well. In practice, we demonstrate that we can reduce AMG solver time by up to $27\%$ with minimal changes to existing PolyDG and VEM software libraries.
\end{abstract}

\maketitle



\section{Introduction}
Multigrid methods are state-of-the-art iterative solvers for large, sparse, linear systems stemming from Finite Element discretizations of Partial Differential Equations (PDEs)  \cite{trottenberg2000multigrid}. A hierarchy of ``coarser" problems is a key ingredient of Multigrid methods, enabling fast error reduction across multiple levels within iterative Krylov-based algorithms  \cite{saad2003iterative}. Within the framework of multigrid methods, Geometric Multigrid methods  \cite{wesseling2004introduction, bramble2019multigrid} take advantage of the underlying PDE and computational domain. Indeed, the sequence of ``coarser" problems is built based on employing coarser grids 
which allows the method to use geometric information to define interpolation and restriction operators between levels. In contrast, Algebraic Multigrid (AMG) methods  \cite{brandt1984algebraic, brandt1986algebraic} do not require knowledge of the underlying problem. AMG methods operate solely on the linear system of equations, analyzing the matrix structure (interpreting it as a weighted graph) to generate coarse levels and the corresponding transfer operators. Therefore, AMG methods can be used as ``black box" or for problems posed on complicated geometries and/or with highly heterogeneous materials, where computing the sequence of coarser meshes can be challenging.\\

Over the past ten years, there has been significant progress in developing new discretization methods for PDEs, commonly referred to as polytopal methods. Polytopal methods extend the Finite Element approach to support computational domains partitioned into general polygons (in two dimensions) or polyhedra (in three dimensions).  
In this paper we focus on Virtual Element Methods  \cite{beirao2013basic,BBM_2013,AHMAD2013,beirao2014hitchhiker} and polytopal Discontinuous Galerkin methods \cite{Antonietti_Brezzi_Marini_2009, Bassi_et_al_2012, antonietti2013hp,cangiani2014hp} which belong to the family of polytopal methods; we also refer the reader to the review papers \cite{BBM_acta_2023,antonietti2016review} and monographs \cite{Antonietti2022VEM, CangianiDongGeorgoulisHouston_2017}, and the references therein.
Other remarkable examples of polytopal methods include Hybrid High-Order Methods~ \cite{DiPietroErnLemair_2014, di2020hybrid}, Hybridizable Discontinuous Galerkin~ \cite{cockburn2010hybridizable, nguyen2010hybridizable}, Mimetic Finite Differences  \cite{lipnikov2014mimetic,BeiraodaVeiga_Lipnikov_Manzini_2014}, and weak Galerkin methods \cite{wang2013weak, wang2016weak}. 
While the development of polytopal discretizations has undergone rapid growth over the last ten years, the development of efficient solvers for the resulting large, sparse (and often ill-conditioned) linear systems is still in its infancy. A key advantage of polytopal methods is that they naturally support agglomeration, i.e., a coarser mesh of arbitrarily shaped elements is obtained by merging elements of an underlying finer grid. 
Mesh agglomeration lacks a direct equivalent in traditional finite elements and is crucial for reducing computational complexity in domains with ``small" inclusions or complex embedded layers and structures. It also forms the basis of constructing efficient geometric multigrid solvers and adaptive algorithms. Indeed, agglomeration (possibly automatically generated by Graph Neural Networks \cite{antonietti2024agglomeration, antonietti2025magnet} or by R-tree algorithms \cite{Feder2025}) form the basis for automatically generating the sequence of nested/non-nested coarser levels in multilevel solvers. For agglomeration-based geometric multigrid solvers for PolyDG and VEM discretizations, we refer to  \cite{antonietti2017multigrid, pan2022agglomeration} and  \cite{antonietti2018multigrid, antonietti2023agglomeration}, respectively. 
While the flexibility of polytopal methods offers advantages in geometric multigrid solvers, this generality also complicates AMG. Indeed, the sparsity pattern of the matrix may be more irregular; elements quality, face/edge connectivity, etc., may vary more. Therefore, interpolation/prolongation operators, and algebraic aggregation strategies have to handle more ``variability".
Since its introduction \cite{brandt1983algebraic}, AMG methods have been widely studied from both algorithmic and theoretical perspectives. 
Developments such as smoothed aggregation  \cite{vanek1996algebraic, van2001convergence} and adaptive variants  \cite{brezina2006adaptive} have broadened its applicability to many finite element discretizations and physical models. The theoretical foundations of AMG methods can be found in \cite{brandt1983algebraic, ruge1987algebraic,falgout2004generalizing,falgout2005two,zikatanov2008two}; we also refer to the review paper  \cite{xu2017algebraic}. 
It is well known that AMG performance depends on parameter tuning, including coarsening strategies, interpolation operators, and the choice of the smoother. AMG faces further challenges when applied to solve the linear systems of equations stemming from polytopal discretizations. The greater variability makes such a tuning both more difficult and more critical. For example, a poor choice of the strong threshold parameter—governing the coarsening step—can severely degrade solver performance. On the other hand, optimal tuning can dramatically accelerate convergence. Analogously, because AMG relies heavily on smoothers to damp high-frequency errors, the choice of the smoother is also critical.
On the other hand, while parameter tuning is known to be critical for achieving good AMG performance, manual tuning is time-consuming and problem- and discretization-dependent.\\

In this paper, we propose using Artificial Neural Networks (ANN) to automate on-the-fly parameter selection for AMG, thereby significantly improving computational efficiency and avoiding manual tuning. 
Our framework belongs to the broad category of learning-augmented numerical methods, which, alongside hybrid, physics-aware, and hidden-dynamics discovery learning frameworks, is among the pillars of the scientific machine learning for PDEs research paradigm, see, e.g., \cite{Bertozzi2012,BruntonKutz2016,BruntonKutz2018,Bertozzi2019,Raissi2019,Quarteroni2019,BruntonKutz2019,Karniadakis2021,Koumoutsakos2020,Zuazua_2022,Zuazua_2023,BELLOMO2024,Quarteroni2025}.
More in detail, we propose a novel ANN-AMG method for the linear systems of equations stemming from VEM and PolyDG discretizations of second-order elliptic problems.
The main idea is that the deep learning algorithm predicts on the fly the optimal AMG strong threshold parameter driving the algebraic coarsening as well as the optimal choice of the smoother. Conceptually, our method combines the advantages of adaptive AMG  \cite{brezina2006adaptive} with that of calibrated AMG. We use Artificial Neural Networks to automatically tune the strong threshold parameter and the best smoother (from a set of given relaxation methods), thereby minimizing the solver time-to-solution. To achieve our goal, we first interpret the matrix of the algebraic system as a grayscale image.
Next, a convolutional and pooling layer generates a compact, multi-channel representation that preserves key structural features while simultaneously reducing computational cost. The proposed ANN-AMG method is entirely non-intrusive, requiring no modifications to existing PDE solver implementations, nor to AMG solvers, thereby guaranteeing compatibility with existing libraries and preserving parallel AMG implementations  \cite{yang2002boomeramg}. Our approach builds upon earlier work that integrated deep learning with multigrid to accelerate iterative solvers \cite{antonietti2023accelerating, caldana2024deep} for Conforming Finite Element Discretizations. However, it extends these ideas in several directions necessary for handling polytopal discretizations. More specifically, we propose optimizing the AMG performance with respect to both the choices of the threshold parameter and the smoother. We also introduce a novel pre-processing step to ensure robustness of training data. We also propose a novel ANN architecture with an additional output measuring the model's prediction confidence. Finally, we propose an improved pooling strategy to capture better variability in the matrix structure and a new training acceleration strategy based on layer freezing. We validate our approach through extensive numerical experiments on both two- and three-dimensional differential diffusion and elasticity problems, discretized using PolyDG and VEM on polytopal meshes.
Results demonstrate that our ANN-AMG method consistently lowers the computational costs by up to $30\%$ compared to classical AMG with default parameters. We note that VEM and PolyDG discretizations are paradigmatic examples of polytopal methods. In VEM, the degrees of freedom are associated with ``geometric" entities such as vertices, edges, faces, or internal moments. In contrast, in PolyDG methods, the local approximation space is ``geometry-agnostic," since a local modal polynomial expansion is used. Consequently, our results demonstrate that the proposed ANN-AMG algorithm seem to achieve acceleration regardless of whether the discretization space is ``virtual and skeleton-based" or ``plain polynomial and modal".
Furthermore, considering both diffusion and elasticity problems serves to show that ANN-AMG is robust for both scalar and vector-valued equations. We remarks that vector-valued equations are  particularly challenging within the AMG framework due to the need for suitable aggregation strategies. Finally, we emphasize that the computational overhead of the ANN forward pass is negligible compared to the acceleration it provides to the solver. Since the training phase is performed offline, the ANN-AMG algorithm can be generalized to a wide range of differential problems, leading to symmetric and positive definite algebraic systems.\\

The remainder of this work is organized as follows. In Section~\ref{sec:AMG}, wereview the basic principles of AMG methods. In Section~\ref{sec:polytopal methods}, we discuss the PolyDG and VEM discretization for both diffusion and linear elasticity problems. Section~\ref{sec:ANN} introduces our deep ANN-AMG algorithm, detailing its architecture and the matrix-to-image pooling representation. In Section~\ref{sec:numerical-results} we test our ANN-AMG algorithm on a wide set of numerical benchmarks carried out on both two- and three-dimensional diffusion and linear elasticity problems. Finally, Section~\ref{sec:conclusions} concludes with a discussion on the obtained results and outlines directions for future research.

\section{Algebraic Multigrid Methods}\label{sec:AMG}
This section introduces the key steps required to construct the AMG method's hierarchy of grids and operators. We first describe the coarse–fine partitioning strategy that determines the grid structure, and then define the interpolation operator that transfers information between levels.

AMG  \cite{brandt1983algebraic, ruge1987algebraic, xu2017algebraic} is an iterative method for solving large, sparse symmetric positive definite (SPD) linear systems of the form 
\begin{equation}\label{eq:system}
A \mathbf{u} = \mathbf{f},
\end{equation}
where $A \in \mathbb R^{n \times n}$ and $\mathbf{f}\in \mathbb R^n$. 
The method constructs a hierarchy of smaller systems to reduce error across different frequency modes efficiently. Its main components are a smoother, grid transfer operators, and a sequence of coarse-grid operators.

The smoother is a simple iterative solver applied for $\nu$ steps, $\mathbf{u}_{k+1} = \mathbf{u}_{k} + S(\mathbf{f} - A \mathbf{u}_{k}), \; k \geq 0$, starting from an initial guess $\mathbf{u}_0$. The matrix $S \in \mathbb{R}^{n \times n}$ defines the method (smoother) that is designed to damp high-frequency components of the error. This operation is denoted as $\texttt{smooth}^\nu(A, S, \mathbf u_0, \mathbf f)$.

To address low-frequency error, the problem is transferred to a coarser grid. This is achieved using interpolation operators $I_{k}^{k-1} \in \mathbb R^{n_{k-1} \times n_{k}}$ and restriction operators $ I_{k-1}^{k} \in \mathbb R^{n_{k} \times n_{k-1}}$. For the SPD case, these operators are defined recursively for levels $k=1, \dots, M-1$ alongside the coarse-grid system matrices $ A^{(k+1)}$:
\begin{equation}
    {I}_k^{k+1} = ( I_{k+1}^{k})^\top, \quad  A^{(k+1)} =  I_k^{k+1}  A^{(k)}  I_{k+1}^{k}, \quad \forall \, k=1,...,M-1, \qquad  A^{(1)} =  A.
\label{eq:op-def-spd}
\end{equation}
The system sizes decrease at each level, $n = n_1 > n_2 > \dots > n_M$. These components are assembled in a recursive procedure, such as the V-cycle shown in Algorithm~\ref{a:v-cycle}. The V-cycle projects the residual equation onto a coarse space ($ A \mathbf e = \mathbf r = \mathbf f -  A \mathbf u_\nu$), solves the problem there, and interpolates the correction back to the fine space ($\mathbf u_\nu + \mathbf e$). 
It uses pre-smoothers $S_1^{(k)}$ for $\nu_1$ steps and post-smoothers $S_2^{(k)}$ for $\nu_2$ steps at each level $(k)$.

\begin{algorithm}[t]
\caption{One Iteration of the V-cycle of the AMG method \newline $\mathbf{u}^{(k)} = \texttt{vcycle}^{k}(\mathbf{u}^{(k)}, \mathbf{f}^{(k)}, \{({A}^{(j)}, {S}_1^{(j)}, {S}_2^{(j)})\}_{j=k}^M, \{({I}^{j+1}_{j}, {I}^{j}_{j+1})\}_{j=k}^{M-1}, \nu_1, \nu_2)$}
\label{a:v-cycle}
\begin{algorithmic}[1]
    \If{$k = M$}
        \State $\mathbf{u^{(M)}} = \texttt{gaussian\_elimination}({A}^{(M)}, \mathbf f^{(M)})$\;
    \Else
        \State $\mathbf{u}^{(k)} \leftarrow \texttt{smooth}^{\nu_1}({A}^{(k)},  S_1^{(k)}, \mathbf{u}^{(k)}, \mathbf{f}^{(k)})$
        \State $\mathbf{r}^{(k+1)}$ $\leftarrow$ $ I_k^{k+1} (\mathbf{f}^{(k)} - {A}^{(k)} \mathbf{u}^{(k)})$   
        \State $\mathbf{e}^{(k+1)}$ $\leftarrow$ $\texttt{vcycle}^{k+1}($\parbox[t]{.6\linewidth}{${\mathbf u}^{(k)}, \mathbf{f}^{(k)}, $ $ \{({A}^{(j)},{S}_1^{(j)},{S}_2^{(j)})\}_{j=k+1}^M, \{({I}^{j+1}_{j}, {I}^{j}_{j+1})\}_{j=k+1}^{M-1}, \nu_1, \nu_2)$}
        \State $\mathbf{u}^{(k)}$ $\leftarrow$ $\mathbf{u}^{(k)} + {I}_{k+1}^k \mathbf{e}^{(k+1)}$
        \State $\mathbf{u}^{(k)} \leftarrow \texttt{smooth}^{\nu_2}({A}^{(k)},  S_2^{(k)}, \mathbf{u}^{(k)}, \mathbf{f}^{(k)})$
    \EndIf
\end{algorithmic}
\end{algorithm}

\subsection{Coarse-Fine Partitioning}
The core of AMG is the automated construction of the interpolation operator based on the entries of the matrix $ A$. For clarity, we will describe the construction for a two-level system, omitting the level-identifying superscripts and subscripts $k$. First, the set of variables $\mathcal{N} = \{1, ..., n\}$ is partitioned into a set of coarse-grid variables $\mathcal{C}$ and a set of fine-grid variables $\mathcal{F}$. Variables in $\mathcal{C}$ will exist on the next coarser grid, while variables in $\mathcal{F}$ will be interpolated from them.
To perform the $\mathcal{C}/\mathcal{F}$ splitting, one must first quantify the coupling between variables based on the matrix entries. Namely, it is crucial to know when a variable $i$ can be interpolated from a variable $j$, or formally, when the variable $i$ \emph{strongly depends} on the variable $j$.

\begin{definition}
\label{def:strong_conn}
Given a threshold parameter $0 < \theta \leq 1$, the set of variables on which variable $i$ strongly depends, denoted $\mathcal{S}_i$, is
\begin{equation*}
\mathcal{S}_i = \{j \neq i: - a_{ij} \geq \theta \, \max_{l \neq i} \, \{ - a_{il} \}, j=1,...,n_k \}.
\end{equation*}
\end{definition}
Conversely, we define the set of variables that are strongly influenced \emph{by} the variable $i$. This ``transpose" set, $\mathcal{S}_i^\top$, is given by $\mathcal{S}_i^\top = \{j: i \in \mathcal{S}_i, j = 1, ..., n_k\}$. Both sets are fundamental for constructing the coarsening and interpolation operators.

We use the CLJP (Cleary-Luby-Jones-Plassman) algorithm to partition the grid. This algorithm models the strong dependencies as a graph $\mathcal{G}=(\mathcal{N}, \mathcal{E})$, where an edge $(i, j) \in \mathcal{E}$ exists if $j \in \mathcal{S}_i$. It then iteratively selects an independent set of nodes to form $\mathcal{C}$. A weight $\eta(i) = |\mathcal{S}_i^\top| + \tilde \eta$ is assigned to each node $i$, where $\tilde \eta$ is a small random number to break ties. The algorithm proceeds until all nodes are classified as either $\mathcal{C}$ or $\mathcal{F}$. This process is repeated to build the operator hierarchy until the grid size $n_k$ is below a threshold, here set to two.

\subsection{Interpolation Operator}
With the $\mathcal{C}/\mathcal{F}$ partition established, the interpolation operator ${I}_{k+1}^k$ is defined as follows. For any vector $\mathbf{x} \in \mathbb R^{n_{k+1}}$ on the coarse grid, its interpolated counterpart on the fine grid is given by:
\begin{equation}
({I}_{k+1}^k \mathbf{x})_i = 
\left\{
\begin{matrix}
    (\mathbf{x})_i & \text{ if } i \in \mathcal{C}^k, \\
    \sum_{j \in \mathcal{C}_i^k} \omega_{ij}^k(\mathbf{x})_j & \text{ if } i \in \mathcal{F}^k,
\end{matrix}
\right.
\label{eq:interp-def}
\end{equation}
where $\mathcal{C}_i^k = \{j \in \mathcal{C}^k: a_{ij} \neq 0\}$ is the set of coarse neighbors of a fine point $i$, and $\omega_{ij}^k$ are the interpolation weights.

The weights are derived from the assumption that the error $\mathbf{e}$ is smooth, meaning $({A}\mathbf{e})_i \approx 0$ for all $i$. This condition can be written as:
\begin{equation*}
\sum_{j\in \mathcal{C}_i} a_{ij} (\mathbf{e})_j+ \sum_{j\in \mathcal{D}_i^s} a_{ij} (\mathbf{e})_j + \sum_{j\in \mathcal{D}_i^w} a_{ij} (\mathbf{e})_j = 0, \quad \forall i = 1, ..., n,
\end{equation*}
where $\mathcal{D}_i^s = \mathcal{F} \cap \mathcal{S}_i$ are the strongly connected fine-grid neighbors and $\mathcal{D}_i^w = \{j \in \mathcal{N}: a_{ij}\neq 0, j \notin \mathcal{S}_i\}$ are the weakly connected neighbors. This leads to the following formula for the weights:
\begin{equation}
\omega_{ij} = - \frac{1}{a_{ij} + \sum_{l \in \mathcal{D}_i^w}a_{il}} \left(a_{ij}+\sum_{l \in \mathcal{D}_i^s}\frac{a_{il}\hat{a}_{lj}}{\sum_{m\in\mathcal{C}_i}\hat{a}_{lm}}\right)
\label{eq:interp_w_def}
\end{equation}
where $\hat{a}_{ij}$ is defined as $a_{ij}$ if $a_{ij}a_{ii}\leq0$ and zero otherwise. The complete AMG setup procedure is summarized in Algorithm~\ref{a:amg}.

\begin{algorithm}[t]
    \caption{AMG algorithm \newline $\mathbf{u} = \texttt{AMG}(\mathbf{u}, {A}, \mathbf{f}, \theta, \{({S}_1^{(j)}, {S}_2^{(j)})\}_{j=k}^M, \nu_1, \nu_2, N_{max}, tol)$}
    \label{a:amg}
    \begin{algorithmic}[1]
        \State build $\{\mathcal{C}^k, \mathcal{F}^k\}_{k=1}^M$ using $\theta$ by means of CLJP
    \State build the operators $\{{I}^{j}_{j+1}\}_{j=k}^{M-1}$ employing Eq.~(\ref{eq:interp-def}) and Eq.~(\ref{eq:interp_w_def})
        \State build the operators $\{{I}^{j+1}_{j}\}_{j=k}^{M-1}$ and $\{({A}^{(j)}\}_{j=k}^M$ by means of Eq.~(\ref{eq:op-def-spd})
        \While{$k < N_{max} $ \textbf{and} $ \norm{{A}\mathbf{u}-\mathbf{f}}/\norm{\mathbf{f}} < tol$}
            \State{$\mathbf u \leftarrow  \texttt{vcycle}^{1}(\mathbf{u}, \mathbf{f}, \{({A}^{(j)}, {S}_1^{(j)}, {S}_2^{(j)})\}_{j=k}^M, \{({I}^{j+1}_{j}, {I}^{j}_{j+1})\}_{j=k}^{M-1}, \nu_1, \nu_2)$}
            \State $k \leftarrow k + 1$
        \EndWhile
    \end{algorithmic}
\end{algorithm}

\section{Model Problems their Polytopal Discretizations}\label{sec:polytopal methods}

In this section, we introduce the model problems under consideration, the basic notation for PolyDG and VEM discretizations, and discuss the corresponding discrete formulations for diffusion and linear elasticity differential problems that lead to linear systems of the form \eqref{eq:system}.\\

Let $\Omega \subset \mathbb{R}^d$, $d=2,3$ be an open, bounded polygonal/polyhedral domain. The first model problem that we consider is the diffusion equation
\begin{equation}\label{eq:model_problem_diffusion}
\begin{aligned}
-\nabla \cdot (\kappa \nabla u) &= f && \text{in } \Omega,\\
u &= 0 && \text{on } \partial \Omega,\\
\end{aligned}
\end{equation}
where the datum $f\in L^2(\Omega)$ and the function $\kappa : \Omega \to \mathbb{R}$ is bounded and uniformly positive. For the sake of simplicity, in the following we will assume that $\kappa$ is piecewise constant over $\Omega$ with discontinuities aligned with the mesh.  

The second model problem under investigation is the following: find the displacement field $\mathbf{u}: \Omega \to \mathbb{R}^d$ such that
\begin{equation}\label{eq:model_problem_elasticity}
\begin{aligned}
 -\nabla \cdot \btensor{\sigma}(\bvector{u}) &= \bvector{f} \quad &&\text{in } \Omega, \\
\bvector{u} &= \mathbf{0} \quad &&\text{on } \partial \Omega,
\end{aligned}
\end{equation}
where $\bvector{f} : \Omega \to \mathbb{R}^d$ is the body force, and $\btensor{\sigma}(\bvector{u})$ is the Cauchy stress tensor. For isotropic elastic materials, $\btensor{\sigma}(\bvector{u}) = \btensor{\mathbb{C}}:\btensor{\varepsilon}(\mathbf{u})=2\mu\, \btensor{\varepsilon}(\bvector{u}) + \lambda\, \operatorname{tr}(\btensor{\varepsilon}(\bvector{u})) \btensor{I},$ where $\btensor{\varepsilon}(\bvector{u})$ is the strain tensor, $\lambda, \mu > 0$ are the Lamé parameters, and $\btensor{I}$ is the identity matrix. The Lamé parameters $\lambda$ and $\mu$ can be expressed in terms of  the Young's modulus $E$ and Poisson's ratio $\nu$ as follows:
\begin{equation}
\mu = \frac{E}{2(1+\nu)},
\qquad\lambda = \frac{E\nu}{(1+\nu)(1-2\nu)},
\label{eq:lame-young}
\end{equation}
respectively.

We next introduce the basic notation, which is instrumental to our discretizations. Let $\mathcal{T}_h$ be a \emph{polytopal} mesh made of general polygons (for $d=2$) or polyhedra (for $d=3$). We denote such polytopal elements by $K$, define by $h_K$ their diameter, and set $h = \max_{\kappa \in \mathcal{T}_h} h_K$. To deal with polygonal and polyhedral elements, we define an \textit{interface} of $\mathcal{T}_h$ as the intersection of the $(d-1)$-dimensional faces of any two neighboring elements of $\mathcal{T}_h$. If $d=2$, an interface/face is a line segment and the set of all interfaces/faces is denoted by $\mathcal{F}_h$.
When $d=3$, an interface can be a general polygon that we assume could be further decomposed into a set of planar triangles collected in the set  $\mathcal{F}_h$.  
Let $F$ be  an interior face shared by two neighboring elements $K^\pm \in \mathcal{T}_h$, i.e., $F = \partial K^+ \cap \partial K^-$.  For regular enough scalar-, vector-, and tensor--valued functions $v$, $\bvector{w}$, and $\btensor{\tau}$,  we define the average and jump operators as
\begin{equation}
    \label{eq:avg_jump_operators}
    \begin{aligned}
    & \jmp{v} = v^+ \mathbf{n^+} + v^- \mathbf{n^-}, \ 
    && \jmp{\bvector{w}} = \mathbf{\bvector{w}}^+ \otimes \bvector{n^+} + \bvector{w}^- \otimes \bvector{n^-}, \ 
    &&
    \\ 
    & \avg{v} = \frac{v^+ + v^-}{2}, \
    && \avg{\bvector{w}} = \frac{\bvector{w}^+ + \bvector{w}^-}{2}, \ && \avg{\btensor{\tau}} = \frac{\btensor{\tau}^+ + \btensor{\tau}^-}{2}.
    \end{aligned}
\end{equation}
 The notation $(\cdot)^{\pm}$ is used to denote the trace on $F$ taken within the interior of $K^\pm$, $\mathbf{n}^\pm$ is the outer unit normal vector to $\partial \kappa^\pm$, and $\bvector{w} \otimes \bvector{n} = \bvector{w}\bvector{n}^T$. On boundary faces, the definition extends as in   \cite{Arnold_Brezzi_Cockburn_Marini_2001}, i.e.,
$$
 \jmp{v} = v \bvector{n},\ \
\avg{v} = v,\ \
\jmp{\bvector{w}} = \bvector{w} \otimes \bvector{n},\ \
\avg{\bvector{w}} = \bvector{w},\ \
\avg{\btensor{\tau}} = \btensor{\tau}.
$$

\subsection{PolyDG Discretization of Diffusion and Linear Elasticity}
We first introduce the scalar broken (discontinuous) polynomial space that serves as the discretization space for PolyDG methods, as
\begin{equation}\label{eq:DGspace}
V_h^{\DG} = \{ v \in L^2(\Omega) : v|_K \in \mathbb{P}^p(K) \qquad \forall K \in \mathcal{T}_h \},
\end{equation}
where $\mathbb{P}^p(K)$ denotes the space of polynomials of total degree at most $p\geq1$. We also need its vector-valued counterpart defined as $\mathbf{V}_h^{\DG}=[V_h^{\DG}]^d$, $d=2,3$.\\

The PolyDG approximation to the model problem \eqref{eq:model_problem_diffusion} reads as follow: find $u_h \in V_h^{\DG}$ such that
\begin{equation}\label{eq:DG_diffusion}
    \mathcal{A}_{\D}^{\DG}(u_h, v_h) = \sum_{K \in \mathcal{T}_h} \int_K f v_h \qquad \forall v_h \in V_h^{\DG}
\end{equation}
where
\begin{align*}
\mathcal{A}_{\D}^{\DG}(u_h, v_h) = &
 \sum_{K \in \mathcal{T}_h} \int_K \kappa \nabla_h u_h \cdot \nabla_h v_h 
- \sum_{F \in \mathcal{F}_h} \int_F \avg{\kappa \nabla_h u_h}  \cdot \jmp{v_h} \\
&- \sum_{F \in \mathcal{F}_h} \int_F \jmp{u_h}  \cdot \avg{\kappa \nabla_h v_h} 
+ \sum_{F \in \mathcal{F}_h} \int_F \sigma \jmp{u_h}  \cdot\,  \jmp{v_h}.
\end{align*}
Here, $\nabla_h$ denotes the element-wise gradient operator, and the penalty function is defined as 
\begin{equation*}
    \sigma|_{F} = \begin{cases}
        \gamma \left\{ \kappa \frac{p_K^2}{h_K} \right\}_{\texttt{H}}, & F \in \mathcal{F}_h^I \\
         \gamma \kappa \frac{p_K^2}{h_K}, & F \in \mathcal{F}_h^B,
    \end{cases}
\end{equation*}
with $\gamma>0$ denoting a penalty coefficient to be properly chosen (large enough) and $\{\cdot \}_\texttt{H}$ denoting the harmonic average.\\

As for the model problem \eqref{eq:model_problem_elasticity},
we introduce
\begin{align*}
\mathcal{A}_{\E}^{\DG}(\bvector{u}_h,\bvector{v}_h) =& 
\sum_{K \in \mathcal{T}_h} \int_K \btensor{\sigma}(\bvector{u}_h):\btensor{\varepsilon}(\bvector{v}_h)
- \sum_{F \in \mathcal{F}_h} \int_F \avg{\btensor{\sigma}(\bvector{u}_h)}  : \jmp{\bvector{v}_h} \\
&- \sum_{F \in \mathcal{F}_h} \int_F \jmp{\bvector{u}_h}  : \avg{\btensor{\sigma}(\bvector{v}_h)} 
+ \sum_{F \in \mathcal{F}_h} \int_F \eta \jmp{\bvector{u}_h}  :  \jmp{\bvector{v}_h},
\end{align*}
with $\eta$ as in Ref.~ \cite{antonietti2016review} [eq. (9)]. Hence, the PolyDG approximation to the elasticity problem \eqref{eq:model_problem_elasticity} reads as follows: find $\bvector{u}_h \in \bvector{V}_h^{\DG}$ such that 
\begin{equation}\label{eq:DG_elasticity}
    \mathcal{A}_{\E}^{\DG}(\bvector{u}_h, \bvector{v}_h) = \sum_{K \in \mathcal{T}_h} \int_K \bvector{f} : \bvector{v}_h
    \qquad \forall \bvector{v}_h \in \bvector{V}_h^{\DG}.
\end{equation}

\begin{remark}{(PolyDG algebraic form)} We observe that each of the  discrete problems \eqref{eq:DG_diffusion}-\eqref{eq:DG_elasticity} yields a linear system of the form \eqref{eq:system}, where $A$ and $\mathbf{f}$ denote the PolyDG stiffness matrix and the right-hand side corresponding to the chosen bilinear form and functional, respectively, and $\mathbf{u}$ is the corresponding vector of expansion coefficients in the selected basis for the discrete spaces $V_h^{\DG}$ and $\bvector{V}_h^{\DG}$, respectively.  
\end{remark} 

\subsection{VEM Discretization of Diffusion and Linear Elasticity}

Following  \cite{beirao2013basic, BeiraoDaVeiga_Brezzi_Marini_2013, BeiraoDaVeiga_Brezzi_Marini_Russo_2016}, in this section, we present the VEM for the diffusion \eqref{eq:model_problem_diffusion} and linear elasticity \eqref{eq:model_problem_elasticity} problems, focusing for the sake of simplicity to the two-dimensional case; we refer to  \cite{BeiraoDaVeiga_Dassi_Russo_2017} and  \cite{BeiraoDaVeiga_Lovadina_Mora_2015} for the details on the three-dimensional extension.\\

For any $K \in \mathcal{T}_h$, we first introduce the following scalar-  and vector-valued \emph{local} virtual element spaces on $K$ is given by
\begin{align*}
& V_h^{\VEM}(K)=\{v\in H^1(K):\, v|_F\in\mathbb{P}_p (F)\ \forall F\subset\partial K,\ \Delta v\in\mathbb{P}_{p-2}(K)\},
\end{align*}
and
\begin{equation*}
  \bvector{V}_h^{\VEM}(K)=[V_h^{\VEM}(K)]^2,
\end{equation*}
with the convention that for $p=1$, $\mathbb{P}_{-1}=\{0\}$.
From the above definitions, we introduce the global spaces $V_h^{\VEM}$ as $\bvector{V}_h^{\VEM}$ as
\begin{align*}
&V_h^{\VEM}=\{v\in H^1_0(\Omega):\ v|_K\in V_h^{\VEM}(K)\;  \forall K \in \mathcal{T}_h\}, && \bvector{V}_h^{\VEM}=[V_h^{\VEM}]^2.
\end{align*}
A common choice of degrees of freedom for $V_h^{\VEM}$, see, e.g., Ref.~ \cite{BeiraoDaVeiga_Brezzi_Marini_Russo_2014}, is:
\begin{enumerate}
\item values of $v$ at the vertices of $K$,
\item on each $F\subset \partial K$ the moments of $v$ up to order $p-2$, for $p\ge 2$,
\item the internal moments of $v$ up to order $p-2$, for $p\ge2$.
\end{enumerate}
This choice of degrees of freedom is unisolvent for the space $V_h^{\VEM}$ and allows the computation of the following projection operator
\begin{align*}
& \Pgrad{p}:V_h^{\VEM}(K)\to\mathbb{P}_p(K) && \int_K\nabla(\Pgrad{p} v - v)\cdot\nabla q =0 && \forall q\in\mathbb{P}_p(K),\\
\end{align*}
with an additional constraint fixing the constant mode (e.g., element mean). 
For $\bvector{V}_h^{\VEM}(K)$, we can define the degrees of freedom analogously (componentwise) and the elastic energy projector as
\begin{align*}
&\bvector{\Pgrad{p}}:\bvector{V}_h^{\VEM}(K) \to [\mathbb{P}_{p}(K)]^2
&&\int_K \left( \btensor{\sigma}(\bvector{\Pgrad{p}}\bvector{v}-\bvector{v})\right):\boldsymbol{\varepsilon}(\bvector{q}) = 0 
&& \forall\, \mathbf{q}\in [\mathbb{P}_{p}(K)]^2,
\end{align*}
together with the orthogonality/rigid-motion constraint to fix rigid body modes (so that the projector is unique). \\

The Virtual Element discretization of \eqref{eq:model_problem_diffusion} reads:
Find $u_h\in V_h^{\VEM}$ such that
\begin{align*}
\mathcal{A}_{\D}^{\VEM}(u_h,v_h) = \sum_{K\in\mathcal{T}_h} \int_K f_h v_h \quad \forall v_h\in V_h^{\VEM}, 
\end{align*}
where $f_h$ is a suitable polynomial projection of $f$, which is computable using the available degrees of freedom, and where 
\begin{align*}
\mathcal{A}_{\D}^{\VEM}(w_h,v_h) = & \sum_{K \in \mathcal{T}_h} \int_K \kappa_K \nabla(\Pgrad{p} w)\cdot\nabla(\Pgrad{p} v)
+\\ & \sum_{K \in \mathcal{T}_h}  S_{\D}^K\big((I-\Pgrad{p})w,(I-\Pgrad{p})v\big), \qquad
 w,v \in V_h^{\VEM},
\end{align*}
where $S_{\D}^K(\cdot,\cdot)$ is a \emph{computable} symmetric positive definite stabilization form acting on the
``non-polynomial'' part of $(I-\Pgrad{p})$ (typical choices for $S_{\D}^K$ include scaled inner products on the degrees of freedom, the so-called ``dofi-dofi" stabilization).\\

As for the Virtual Element discretization of the linear elasticity problem \eqref{eq:model_problem_elasticity}. 
Find $\bvector{u}_h\in \bvector{V}_h^{\VEM}$ such that
\begin{align*}
\mathcal{A}_{\E}^{\VEM}(\bvector{u}_h,\bvector{v}_h) =  \sum_{K\in\mathcal{T}_h} \int_K \bvector{f}_h \cdot \bvector{v}_h|_K
\quad \forall \bvector{v}_h \in \bvector{V}_h^{\VEM}
\end{align*}
where $\bvector{f}_h$ is a suitable polynomial projection of $f$, which is computable using the available degrees of freedom.
The computable bilinear form $\mathcal{A}_{\E}^{\VEM}(\cdot,\cdot)$ is defined as
\begin{align*}
&\mathcal{A}_{\E}^{\VEM}(\bvector{w}_h,\bvector{v}_h)=\sum_{K \in \mathcal{T}_h}\int_{K} \btensor{\sigma}(\bvector{\Pgrad{p}}\bvector{w}):\btensor{\varepsilon}(\bvector{\Pgrad{p}}\bvector{v})
+ \sum_{K \in \mathcal{T}_h} S_{\E}^K\big((\bvector{I}-\bvector{\Pgrad{p}})\bvector{w},(\bvector{I}-\bvector{\Pgrad{p}})\bvector{v}\big)
\end{align*}
for all $\bvector{w}_h, \bvector{v}_h \in \bvector{V}_h^{\VEM}$
where $S_{\E}^K(\cdot,\cdot)$ is a symmetric positive definite stabilization acting on the kernel of the projector (i.e., the non-polynomial part).

\begin{remark}{(VEM algebraic form)} As before, the Virtual Element discretization of both the diffusion and linear elasticity problem yields a linear system of the form \eqref{eq:system}, where $A$ and $\mathbf{f}$ denote the VEM stiffness matrix and the right-hand side corresponding to the chosen bilinear form and right hand side, respectively, and $\mathbf{u}$ is the corresponding expansion coefficient vectors for the discrete Virtual  Element spaces $V_h^{\VEM}$ or $\bvector{V}_h^{\VEM}$, depending on the discretized differential problem.  
\end{remark} 

\section{The Neural Network Architecture}\label{sec:ANN}
In this section, we provide a detailed description of the neural network architecture and the algorithm used to tune the AMG method automatically. First, we provide an overview of the pipeline, introducing the high-level components and procedures that enable the automatic choice of parameters. Then, in each subsection, we will go into the details and precisely define all the components of our method.

The first choice concerns the parameters to be tuned and the employed metric (that is, a scalar index that establishes  when a certain combination of parameters is better than another).
In Section~\ref{sec:AMG}, we have shown that the strong threshold parameter $\theta$  critically influences the construction of the interpolation operator and thus the whole hierarchy of levels that stand at the basis of each AMG application. Another key choice that heavily influences the behavior of the AMG methods is the smoother. The importance of these parameters is confirmed by empirical evidence when tested for our applications. For these reasons, we aim to tune $\theta$ and the  smoother to minimize the computational cost of the AMG method. We remark that the tuning process can be done by optimizing the choice of other parameters, such as the number of pre- and post-smoothing iterations applied at each level, or the choice of $\mathcal{C}/\mathcal{F}$ splitting algorithm.  

Concerning the performance metric, we ideally want to minimize the total elapsed time needed to solve the linear system. However, this metric has at least two significant drawbacks: it is machine-dependent and it suffers from measurement error. Techniques to tackle these problems are detailed in Section~\ref{sec:data-cleaning}. However, in some instances, statistical analysis shows a strong correlation between the elapsed time and the approximate convergence factor 
\begin{equation}
    \rho = \bigg(\frac{||\mathbf{r}^{(k)}||_2 }{||\mathbf{r}^{(0)}||_2}\bigg)^{N_{it}}
\end{equation}
which measures how quickly the iterative solution contracts towards the exact one. 
Hence, when possible, we prefer to use $\rho$ as the index of the performance. For each test case, we also considered multivariate polynomial models that take into account the size and number of non-zero elements of
$A^{(k)}$ at each level $k$. Unfortunately, none of them  showed statistical significance in predicting the elapsed time $t$.

The optimization step  is performed at each application of the AMG method, that is, every time we solve the system $A \mathbf u = \mathbf f$.  In view of this, it seems quite natural that only the matrix $A$ is employed in the tuning procedure. However, the matrix $A$ is a large, sparse matrix of variable size. Neural networks are not well-equipped to handle this kind of data. For this reason, we use a special kind of pooling, first introduced in Ref.~ \cite{antonietti2023accelerating}, to prune and compress $A$ into a small multi-channel image $\mathbf V \in \mathbb{R}^{m\times m \times f}$. This procedure is outlined, together with the meaning of the hyperparameters $m$ and $f$, in Section~\ref{sec:pooling}.

Hence, our goal is to use a neural network $\mathscr{F}$ to predict the computational cost of solving a system $A \mathbf u= \mathbf f$, where the input of $\mathscr{F}$ is the matrix $A$ and the parameters that we want to tune (the scalar $\theta$ and the categorical variable for the kind of smoother encoded as a one-hot). Although including $\mathbf f$ (or a pooled form of $\mathbf f$) as ANN input might improve per-problem performance, we omit it because AMG constructs its multilevel hierarchy and interpolation operators from $A$ alone. By limiting the network’s input to information derived from $A$, we aim to learn features that generalize across problems rather than overfit to particular right-hand sides; we therefore hypothesize improved generalization. The optimal choice of parameters is then found by solving two optimization problems. The first one is the offline training of the ANN, which enables us to learn an approximate map of the computational cost depending on the problem we are solving ($A$) and the choice of the AMG's parameters. This step is expensive but it is done only once. The second one is an online optimization step that allows us to find the optimal choice of parameters, namely:

\begin{equation}
    \min_{\text{AMG parameters}} \mathscr{F}(A, \cdot \,).
    \label{eq:abstact-online-min}
\end{equation}

This optimization step takes place in the online phase and must be solved each time we apply the algorithm.

While this choice might seem strange at first, we would like to highlight a few advantages. This approach is more data-efficient. Specifically, if the ANN were trained to directly predict the optimal combination of AMG parameters, we would need to solve a  more expensive offline optimization problem, which in turn would require solving many more linear systems $A \mathbf u = \mathbf f$. This procedure would be prohibitively expensive. By contrast, if the ANN is trained to predict the computational cost, then every computational cost measurement taken by solving $A \mathbf u = \mathbf f$ -- for any choice of AMG's parameters -- can be used as a training sample. In this way, we can build a large dataset by solving many inexpensive, small-scale problems, and then supplement it with only a few samples from larger, more costly problems, relying on the ANN's generalization capability.  In other words, by predicting the computational cost instead of the optimal AMG parameters, we build a surrogate model $(A, \text{AMG parameters}) \mapsto c$ which enables a more efficient solution of Eq.~\ref{eq:abstact-online-min} than directly measuring $c$ by solving $A \mathbf u = \mathbf f$. This feature is incredibly useful when the solution of even one linear system may be costly. Moreover, we have empirically found this approach to be more stable with respect to directly learning the optimal value of parameters. Finally, the search space is practically one-dimensional (since the smoother choice is discrete), thus the cost of performing the optimization is small.

\subsection{Handling Measurement Uncertainty}\label{sec:data-cleaning}

In this section, we describe the procedure used to collect the data for training the neural network when the execution time is taken as  performance metric. Since measurements of the elapsed time $t$ are subject to noise, we repeat each measurement multiple times.

To reduce data collection costs while mitigating measurement errors, we adopt the following strategy: the measurement is repeated $r$ times, where $r$ ranges from 2 to 100 and is chosen inversely proportional to the mean elapsed time of the first two measurements. The rationale is that measurement variability arises primarily from operating-system tasks running concurrently with the AMG solver, perturbing CPU load. For larger values of $t$, these fluctuations tend to average out, resulting in lower variance. Empirical evidence supports this assumption, as we observed that the sample variance of repeated measurements decreases inversely proportionally to the elapsed time (for fixed $r$). The reported elapsed time $t$ is then taken as the mean of the $r$ repetitions.

To further improve data quality, we apply a Savitzky-Golay filter  \cite{savitzky1964smoothing}. We employ a window size of 21 and a polynomial of degree 7 for uniformly sampled values of $\theta$. These parameters were determined through manual tuning, by testing multiple combinations and assessing their performance on representative subsets of the dataset. Selection was guided by visual inspection, balancing smoothness and fidelity to the raw data. Importantly, we verified that the filter preserves the positivity of the data and performed manual checks whenever the difference between the minima of the raw and smoothed signals exceeded a prescribed threshold. In such cases, we adjusted the polynomial degree or window size as appropriate. We note that smoothing can alter the position of minima in sharp valleys, as it tends to attenuate high-frequency features. Improper tuning of the filter may significantly degrade predictive accuracy, whereas carefully calibrated smoothing provides accuracy improvements of several percentage points. Nevertheless, in the absence of proper tuning, we observed that omitting smoothing remains preferable to applying it blindly.

All experiments reported here were carried out in a serial manner in a controlled environment, ensuring that no other processes contributed significant CPU load. To enhance training robustness, we normalize execution times to the interval $[0, 1]$. 

Finally, we remark that parallel execution introduces another external factor -- the number of CPU cores -- which can be treated as an additional AMG parameter. However, this extension is left for future work.

\subsection{The Pooling Operator}\label{sec:pooling}
One of the key steps of our algorithm is passing the information contained in $A$ to the neural network. Even if the structure of a sparse matrix is more akin to a graph than a dense matrix, we prefer to use a neural network with structured input (CNN) due to the size of $A$. Namely, to make this process scalable to cases where $A$ has millions of entries, we prefer not to directly apply GNNs to $A$, but instead, we would rather employ a process that can prune and extract information from $A$ much more quickly. In particular, we use a variation of the pooling technique used in CNNs, first introduced in Ref.~ \cite{antonietti2023accelerating}. 

We denote by $\mathbf V = \texttt{pooling}(A, m) \in \mathbb R^{m \times m \times f}$ the pooled representation of $A$, where $m \in \mathbb N$ is a hyperparameter that controls the tensor size and $f$ is the number of features extracted, in our case $f=4$. Letting $q = \lceil n_1 / m \rceil$, the pooled features are defined as
\begin{align*}
& v_{ij1} = \max_{i',j'=1,\ldots,q}{\max \, \{ 0, \tilde a_{iq+i',jq+j'}  \}}, \qquad
& v_{ij2} = \max_{i',j'=1,\ldots,q}{\max \, \{ 0, -\tilde a_{iq+i',jq+j'} \}}, \\
& v_{ij3} = \sum_{i'=1}^{q}\sum_{j'=1}^{q} \tilde a_{iq+i',jq+j'}, \qquad 
& v_{ij4} = \sum_{i'=1}^{q}\sum_{j'=1}^{q} \chi_{(0,+\infty)}(\tilde a_{iq+i',jq+j'}),
\end{align*}
where $\chi$ is the indicator function and $\tilde a_{ij} = a_{ij}\chi_{i\leq n_1,,j\leq n_1}$.

A reference implementation is reported in Algorithm~\ref{a:pooling}. While the pseudocode assumes that $A$ is stored in coordinate (COO) format, the procedure extends to other sparse storage formats. Unlike conventional CNN pooling, which typically extracts a single maximum, here we compute four complementary features within each neighborhood: maximum of positive entries, maximum of negative entries, sum of all entries, and the number of nonzeros. This choice is motivated by the role of positive, negative, and aggregate values in defining the interpolation weights (Eq.~\ref{eq:interp_w_def}), while the nonzero count provides a measure of local sparsity. Although additional features could be incorporated, this selection offers a favorable balance between expressiveness and efficiency. Notably, the features preserve the sparsity structure: pooling a block of zeros yields zero.

The algorithm scales linearly with the number of nonzero entries, i.e., $O(\textit{nnz})$, which for finite element discretizations translates to $O(p\cdot n)$. Empirical tests confirm that its runtime is negligible compared to that of solving the linear system, a prerequisite for its practical relevance. Furthermore, Algorithm~\ref{a:pooling} can be parallelized with minimal effort.

\begin{algorithm}[t]
\caption{Pooling algorithm $\mathbf V=\texttt{pooling}({A},m)$}
\label{a:pooling}
\begin{algorithmic}[1]
\State Access ${A}$ in COO format and extract: $n_1$, \texttt{val}, \texttt{row}, \texttt{col}
\State Initialize $\mathbf{V}$ as an $m \times m \times 4$ dense tensor with zero entries
\State $q \leftarrow n_1/m$, $p \leftarrow n_1 \bmod m$, $t \leftarrow (q+1)p$
\For{$k = 0$ \textbf{to} $\texttt{val.size()}-1$}
\State $i \leftarrow \texttt{row}[k]/(q+1)$ if $\texttt{row}[k]<t$, else $(\texttt{row}[k]-t)/q+p$
\State $j \leftarrow \texttt{col}[k]/(q+1)$ if $\texttt{col}[k]<t$, else $(\texttt{col}[k]-t)/q+p$
\State $v_{ij1} \leftarrow \max{ \max(0, \texttt{val}[k]), v_{ij1}}$
\State $v_{ij2} \leftarrow \max{ \max(0, -\texttt{val}[k]), v_{ij2}}$
\State $v_{ij3} \leftarrow v_{ij3} + \texttt{val}[k]$
\State $v_{ij4} \leftarrow v_{ij4} + 1$
\EndFor
\end{algorithmic}
\end{algorithm}

\subsubsection{Normalization}\label{sec:normalization}

To ensure stable and fast convergence, the tensor $\mathbf V$ is normalized following a logarithmic normalization scheme  \cite{antonietti2023accelerating}. It has been shown that this normalization outperforms standard linear scaling in this context. Namley, for each feature $f$, the transformation reads
\begin{equation}
\hat v_{ijf} = \frac{\log(|v_{ijf}|+1)}{\max_{i,j} |\log(|v_{ijf}|+1)|}\frac{v_{ijf}}{|v_{ijf}|}.
\label{eq:normalization}
\end{equation}
A key property of this normalization is that zeros map to zero, thereby preserving the sparsity pattern of $\mathbf V$.

\subsection{The AMG-ANN Algorithm}\label{sec:amg-ann-algo}
The core of the algorithm is the neural network $\mathscr{F}$ that predicts the computational cost of the AMG method. As mentioned before, the computational cost for us can either be the approximate convergence factor $\rho$ or the normalized and smoothed wall clock time $\bar{t}$.

The architecture of $\mathscr{F}$ is made up of two main components in series. The first is a CNN that analyzes and encodes a flat latent representation for the pooling tensor $\mathbf V$. The hyperparameters that we consider for this component are the size $m$ of $\mathbf V$, the number of convolution blocks (ending with a max-pooling layer), the number of convolutions in each block, and the number and size of filters of each convolution. The activation function we use is the ReLU. 

The latent representation of $\mathbf V$ is concatenated with the parameters of the AMG method we want to tune (the scalar $\theta$ and the kind of smoother encoded as a one-hot of size four), the logarithm of the size of $A$ ($\log(n)$) -- since during the pooling this information is lost -- and, the polynomial degree $p$ of the basis function used in the FE discretization. This information is not strictly needed; indeed, our experiments show that even omitting the parameter $p$, we can obtain a neural network with a similar loss. However, when adding the polynomial degree to the inputs, the training is usually less expensive and requires a shorter phase of hyperparameters' tuning. Hence, our algorithm is not
limited to problems stemming from FE discretizations.

\begin{figure}[!htbp]
    \makebox[\textwidth][c]{\includegraphics[width=0.6\textwidth]{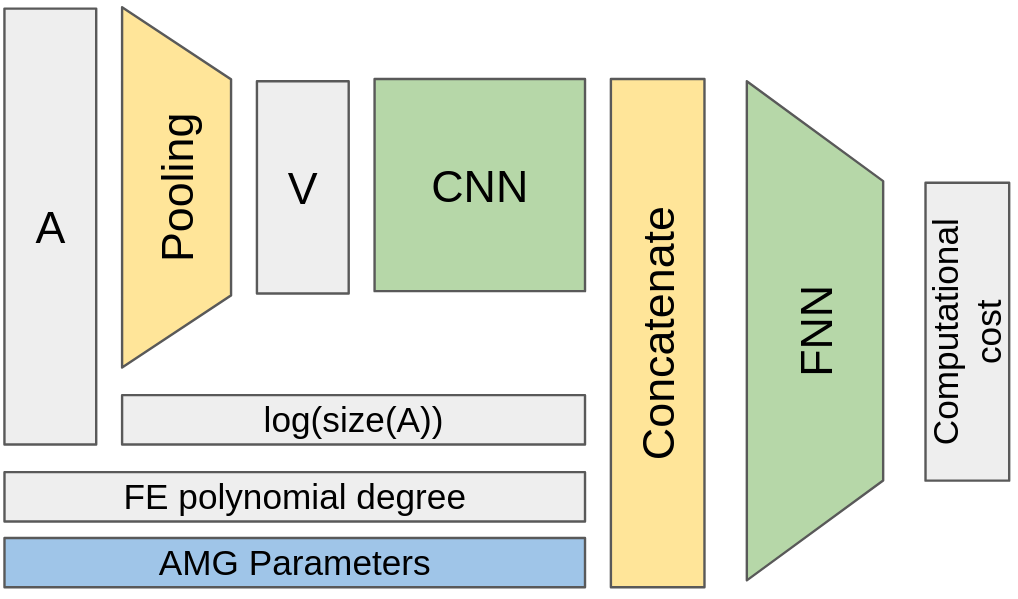}}
    \caption{Architecture of the proposed ANN used to predict the optimal value of the AMG parameters.}
    \label{fig:architecture}
\end{figure}

The second component of $\mathscr{F}$ is a dense feed-forward neural network that takes as input the aforementioned concatenation layer and outputs a prediction of the computational cost. On the last layer, we add a clipping activation function to ensure that the output lies in $[0, 1]$. The hyperparameters of this component are the depth and the width of the network. To sum up, we have
\begin{equation}
    \mathscr{F}(\hat{\mathbf{V}}, \log{n}, p, \theta, b; \boldsymbol{\alpha}) = \tilde{c}
\end{equation}
where 
\begin{equation*}
b \in \mathcal{O}_4 = \left\{ e_i \in \{0,1\}^4 : \sum_{j=1}^{4} (e_i)_j = 1 \right\}
\end{equation*}
is the one-hot representation of the kind of smoother. Namely, it is a vector of four components where only one is not zero. Thus, if the non-zero component is in the $i$-th position, we choose as smoother the $i$-th in the following list: Jacobi with successive over-relaxation (SOR), $\ell^1$-Jacobi  \cite{baker2011multigrid}, $\ell^1$-Jacobi with SOR and Jacobi with FCF relaxation  \cite{hessenthaler2020multilevel}. Moreover, we denote by  $\boldsymbol{\alpha}$ the parameters of the neural network, while  $\tilde{c}$ denotes the approximation of the computational cost $c$. The architecture of the ANN is shown in Figure~\ref{fig:architecture}.
 
Like most deep learning algorithms, our algorithm works in two phases.

\paragraph{Offline phase} We collect, smooth, and normalize the data following the procedure outlined in Section \ref{sec:data-cleaning}. Then we train the neural network with the collected samples. Namely, the input-target training pairs for the supervised training are the couples $((\hat{\mathbf{V}}^{(i)}, \log{n}^{(i)}, p^{(i)}, \theta^{(i)}, b^{(i)}), c^{(i)})$, where the superscript $(i)$ indicates that it is the $i$-th sample of the dataset. In other words, the ANN receives as input the matrix of the system, the threshold parameters, and the choice of the smoothers, and it has as its target the computational cost. As usual, The optimization step  minimizes the MSE error between the predicted cost and the target cost $c$.

The hyperparameters of the neural network are chosen via a Bayesian optimization algorithm  \cite{akiba2019optuna}. The procedure uses a standard 60-20-20 split into train, validation, and test datasets. The data are partitioned at the problem level, meaning that each matrix 
$A$ appears in exactly one of the three datasets. Consequently, when the algorithm is evaluated on the validation or test set, it is required to make predictions on entirely unseen problems. Training has been performed using the AdamW  \cite{loshchilov2017decoupled} optimizer. The initial learning rate and the minibatch size are also tuned hyperparameters. The learning rate is managed via a suitable learning rate schedule to improve convergence speed and stability. Namely, we employ a reduced on-plateau learning rate schedule that halves the learning rate with a manually tuned patience. More details on how data is practically generated are given in Sections~\ref{sec:numerical-results-diffusion},~\ref{sec:numerical-results-elasticity} since it depends on the physics and the FE discretization of the problem.

\paragraph{Online phase} Given as input the matrix $A$ of the linear system to be solved and, optionally, the polynomial degree $p$ of the discretization.
\begin{itemize}
    \item[1] Compress $ A$ by applying the pooling to obtain $\mathbf V$ (Section~\ref{sec:pooling}).
    \item[2] Normalize $ V$ to obtain $\hat{{V}}$ (Section~\ref{sec:normalization})
    \item[3] Obtain $(\theta^*, b^*)$ solving the continuous optimization problem
    \begin{equation}\label{eq:A-theta-map}
        (\theta^*, b^*) = \argmin_{(\theta, b) \in [0, 1] \times \mathcal{O}_4} \mathscr{F}(\hat{\mathbf{V}}, \log{n}, p, \theta, b)
    \end{equation}
    by doing a complete search of the space using the discretization 
    \begin{equation*}
        (\theta^*, b^*) = \argmin_{(\theta, b) \in \texttt{linespace}(0, 1, 101) \times \mathcal{O}_4} \mathscr{F}(\hat{\mathbf{V}}, \log{n}, p, \theta, b)
    \end{equation*}
    \item[4] Use $(\theta^*, b^*)$ as parameters of the AMG method.
\end{itemize}

\subsection{Evaluating the Model}
While obtaining a small MSE test loss between the predicted and target computational cost $c$ during the offline phase is a good indicator that the neural network is learning, it does not measure the actual reduction of cost that our algorithm has on the AMG method.
We remark that the choice of parameters $(\theta, b)$ employed in the AMG to solve the system is subordinate to the map $A \mapsto (\theta^*, b^*)$ defined by the ``Online phase'' algorithm of Section~\ref{sec:amg-ann-algo} and, in particular, to Eq.~(\ref{eq:A-theta-map}).  Hence, we introduce the following quantities of interest. Let $A$ be fixed, and let:

\begin{itemize}
    \item $t_{\textnormal{ANN}}$ be the computational time of the AMG-ANN algorithm, that is by using $\theta=\theta^*$ and $i$-th smoother among SOR-Jacobi, $\ell^1$-Jacobi, SOR-$\ell^1$-Jacobi and FCF-Jacobi, where $i$ is the position of the only non-zero bit of $b^*$; 
    \item $t_{0}$ be the computational time of the AMG method with the default parameter ($\theta=0.25$ in 2D, $\theta=0.5$ in 3D, and SOR-Jacobi smoother);
    \item $t_{\textnormal{MIN}}$ be the computational time of the AMG method with 
\vspace{-3mm}
\begin{equation*}
 (\theta^*, b^*)=\argmin_{ (\theta, b) \in \textnormal{dataset for A}} t(\theta, b; A);
\end{equation*}
\vspace{-4mm}
\item $P = 1 - \frac{t_\textnormal{ANN}}{ t_{0}}$ be the performance index of the AMG-ANN algorithm;
\item $P_\textnormal{MAX} = 1 - \frac{t_{MIN}}{ t_{0}}$ be the best performance of the AMG-ANN algorithm.
\end{itemize}

Moreover, we can compound the quantities over different $A$ and define $P_B$ as the percentage of cases where $P \geq 0$, $P_m$ as the average of $P$, and $P_M$ as the median of $P$. The ratio $P_r = P/P_{MAX}$ gives the accuracy of the neural network, namely a measure of how close the ANN is to performing to the theoretical maximum.

Naturally, all the metrics defined above can be computed also if we employ as a measure of the computational cost the approximate convergence factor $\rho$. Indeed, we can compute $(\theta^*, b^*)$ as explained in the previous section and then compute all the performances indexes ($P$, $P_{MAX}$, etc.) by taking measurements of the wall clock time needed to solve the linear system for each choice of $(\theta, b)$ present in the dataset.
Finally, let us remark that it can happen that, when the matrix $A$ is ill-conditioned, the quantity $t_{0}$ is not well defined since the default values of the parameters of the AMG do not allow the solver to converge. In these cases we set $t_{0}=\infty$ and thus $P=1$. In a few of the following test cases, this is not a rare event. Thus, for the sake of completeness, we also introduce the following additional  performance indexes.
\begin{itemize}
    \item $P_w$ as the percentage of cases where $t_{0}=\infty$;
    \item $t_{0,i}$ as computational time of the AMG method with the default value of $\theta$ and the $i$-th smoother; 
    \item $P^i = 1 - \frac{t_\textnormal{ANN}}{ t_{0,i}}$ the performance with respect to the $i$-th smoother.
\end{itemize}

\section{Numerical Results}\label{sec:numerical-results}
This section presents the numerical experiments designed to evaluate the performance of the proposed AMG-ANN algorithm across a variety of discretizations, problem settings, and dimensions. We first describe the computational framework and solver configuration adopted in our study. This provides the necessary context for interpreting the numerical results presented later. The experiments are then organized into two main categories: the diffusion problem and the linear elasticity problem. Within each category, we report several test cases of increasing complexity, spanning both two- and three-dimensional domains, multiple mesh families, and different discretization schemes (VEM and PolyDG).

A key advantage of our algorithm is its ability to integrate seamlessly with existing software libraries. In particular, we assemble the linear systems using the \texttt{VEM++}  \cite{dassi2025vem++} and \texttt{Vulpes} libraries (\texttt{https://vulpeslib.github.io}). We employ BoomerAMG  \cite{yang2002boomeramg} from the HYPRE library  \cite{falgout2002hypre} as the AMG solver. BoomerAMG is used strictly as a black-box solver; no modifications are introduced except for setting the strong threshold parameter $\theta$ and the smoother. All other parameters remain at their default values. We use AMG as a preconditioner for the Conjugate Gradient Krylov solver of PETSc  \cite{balay2019petsc}.
The stopping criterion is defined by enforcing a relative residual tolerance of $10^{-8}$. \\

Finally, we describe the method used to choose or sample the threshold parameter $\theta$, which is common across all test cases. For small problems (total number of degrees of freedom $n < 20{,}000$), we use a uniform discretization of the interval $[0.05, 0.95]$ with a step size of $0.025$. For moderate problems (total number of degrees of freedom up to $n = 100{,}000$), we employ a coarser discretization with a step size of 0.05. For larger problems (total number of degrees of freedom up to $n > 100{,}000$), we uniformly randomly sample the interval at $10$ or fewer points, depending on the problem size.

\subsection{Numerical Results: Diffusion Problem}\label{sec:numerical-results-diffusion}
\subsubsection{Test Case 1: VEM Discretization of the Diffusion Problem in 2D}
\begin{figure}[!htbp]
    \centering
    \subfloat[Triangular Mesh]{
\includegraphics[width=0.23\textwidth]{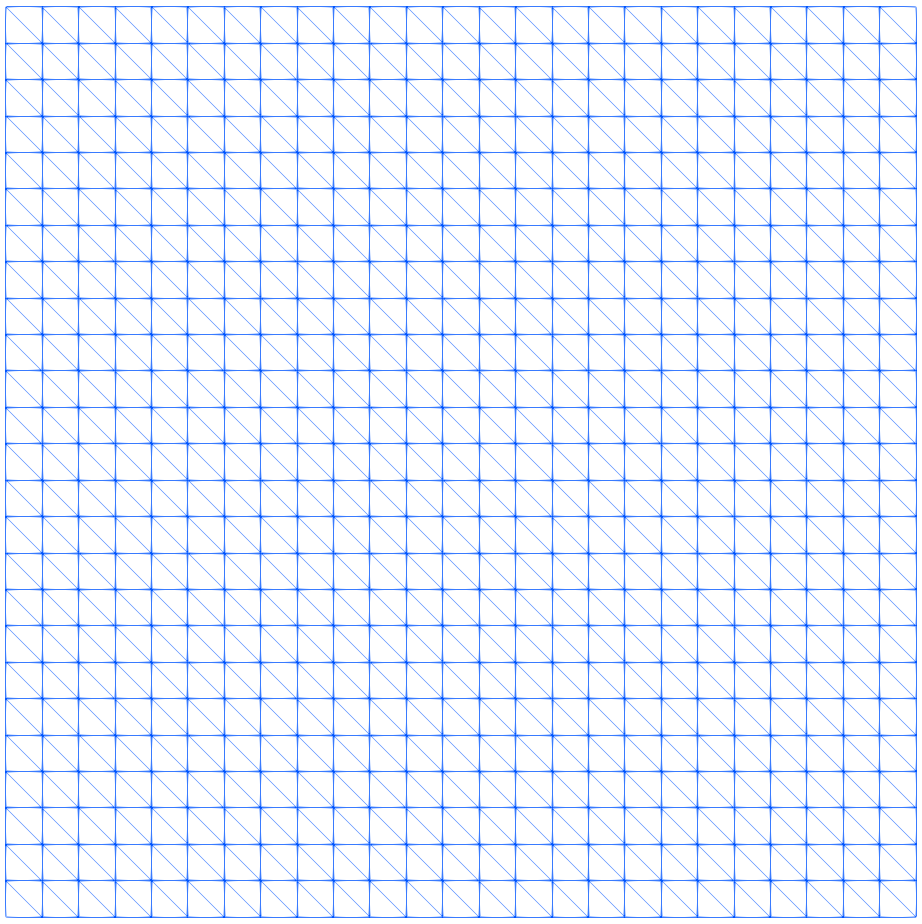}
        \label{fig:tri-mesh}
    }
    \subfloat[Quadrilateral Mesh]{        \includegraphics[width=0.23\textwidth]{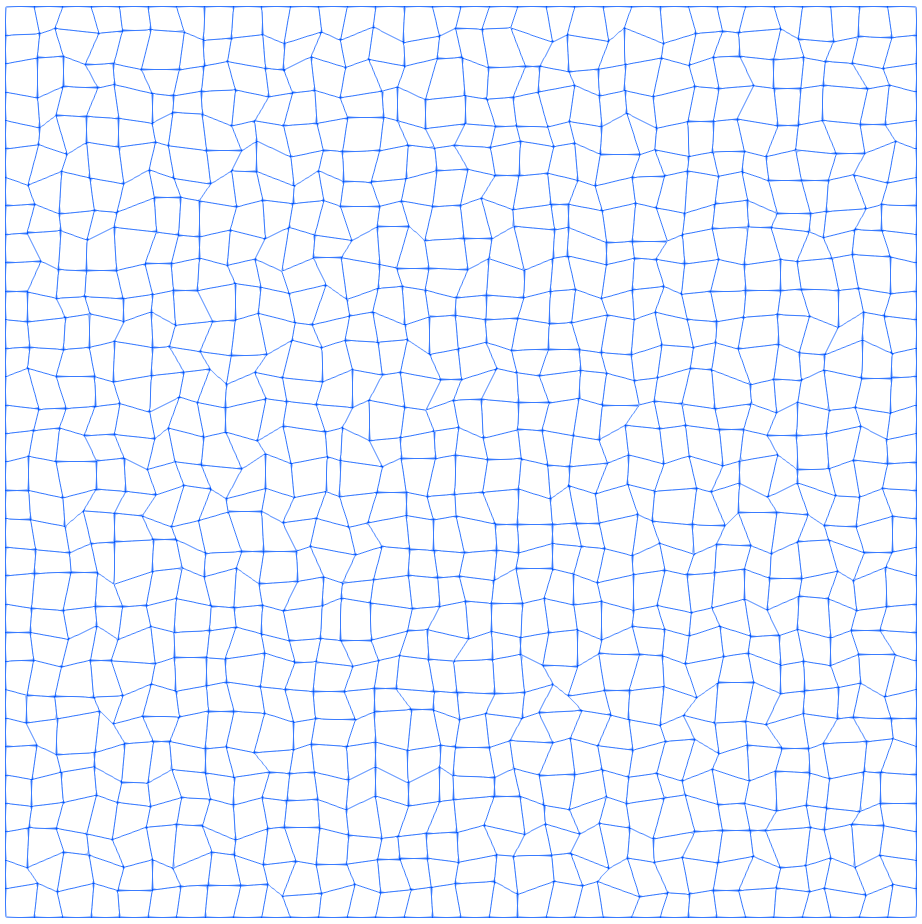}
        \label{fig:quad-mesh}
    }
    \subfloat[Voronoi Mesh]{
        \includegraphics[width=0.23\textwidth]{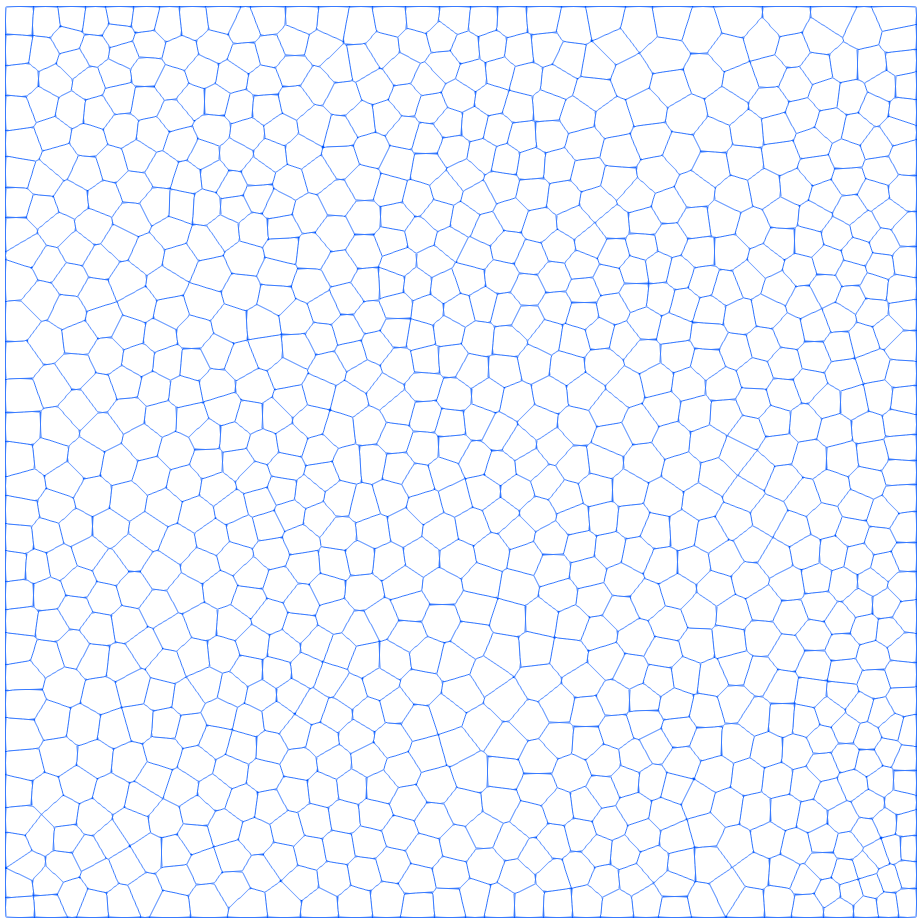}
        \label{fig:voro-mesh}
    }
    \subfloat[Hexagonal Mesh]{
        \includegraphics[width=0.23\textwidth]{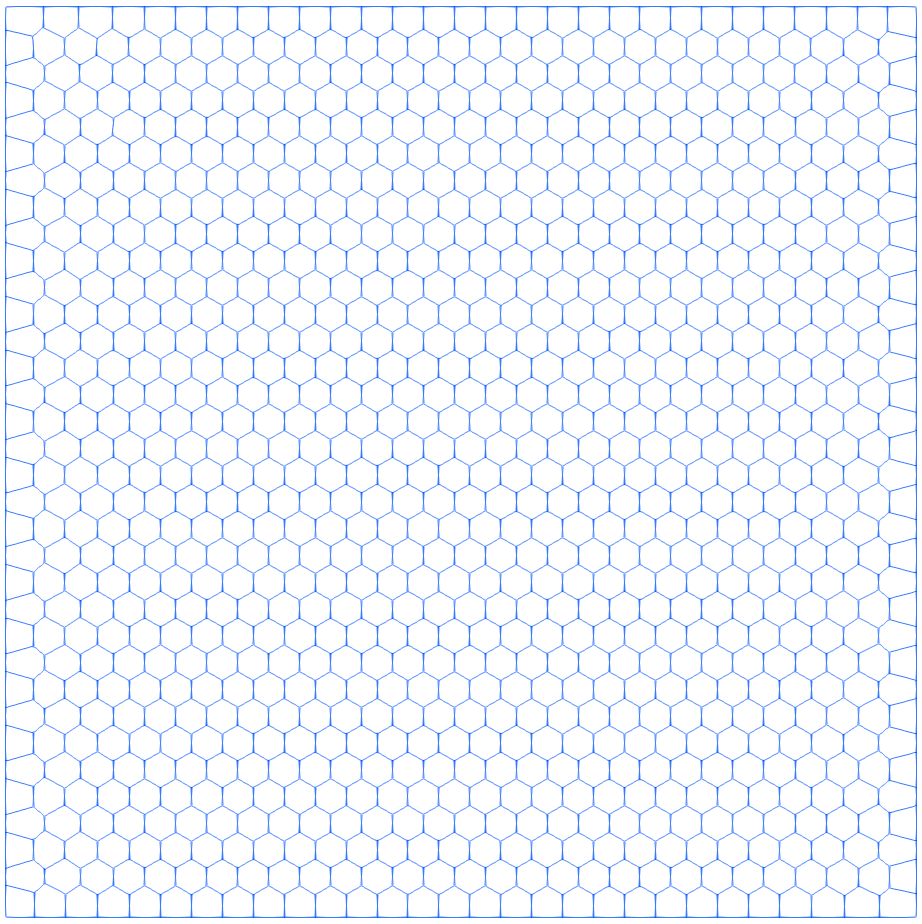}
        \label{fig:hex-mesh}
    }   
    \vspace{1.5em} 
    \subfloat[Triangular\\Sparsity]{
        \includegraphics[width=0.23\textwidth]{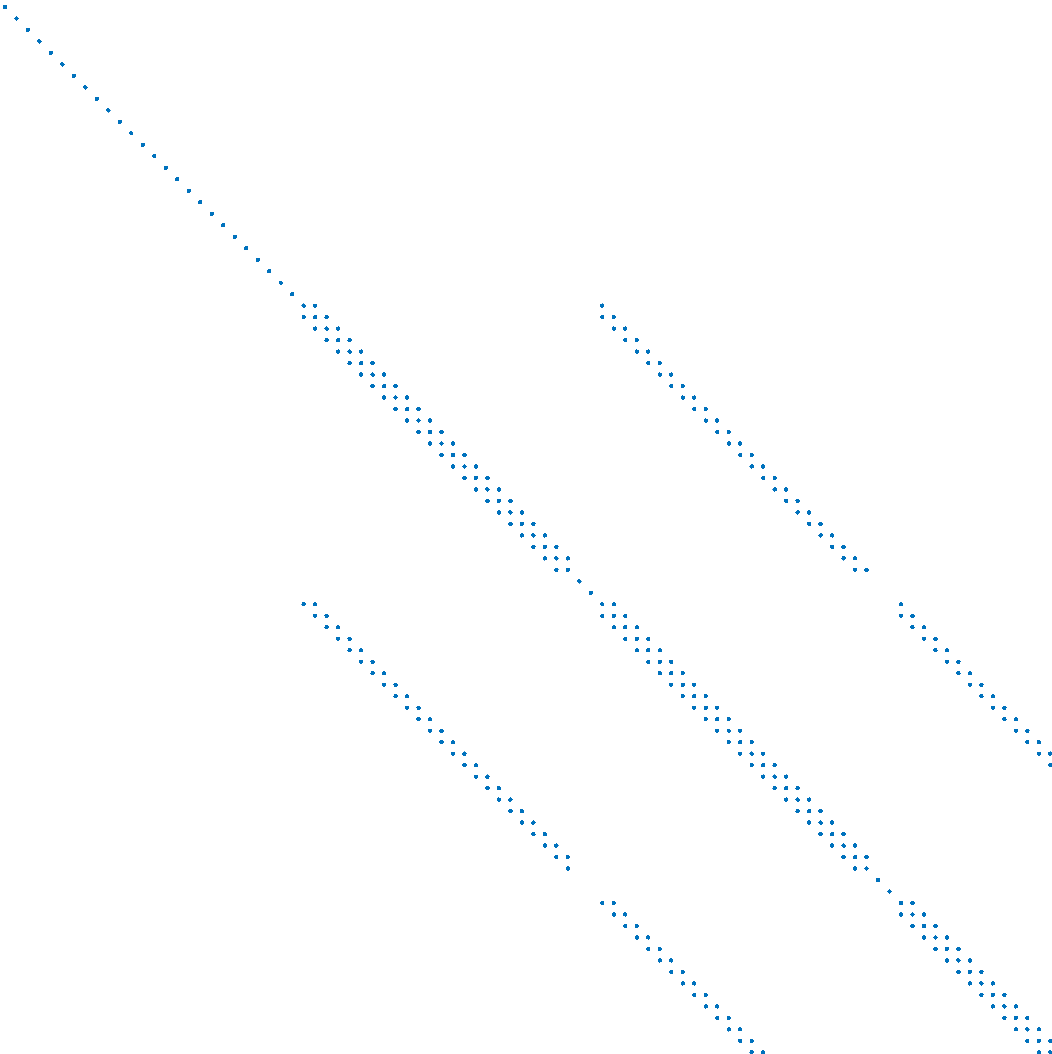}
        \label{fig:tri-sparsity}
    }
    \subfloat[Quadrilateral\\ Sparsity]{
        \includegraphics[width=0.23\textwidth]{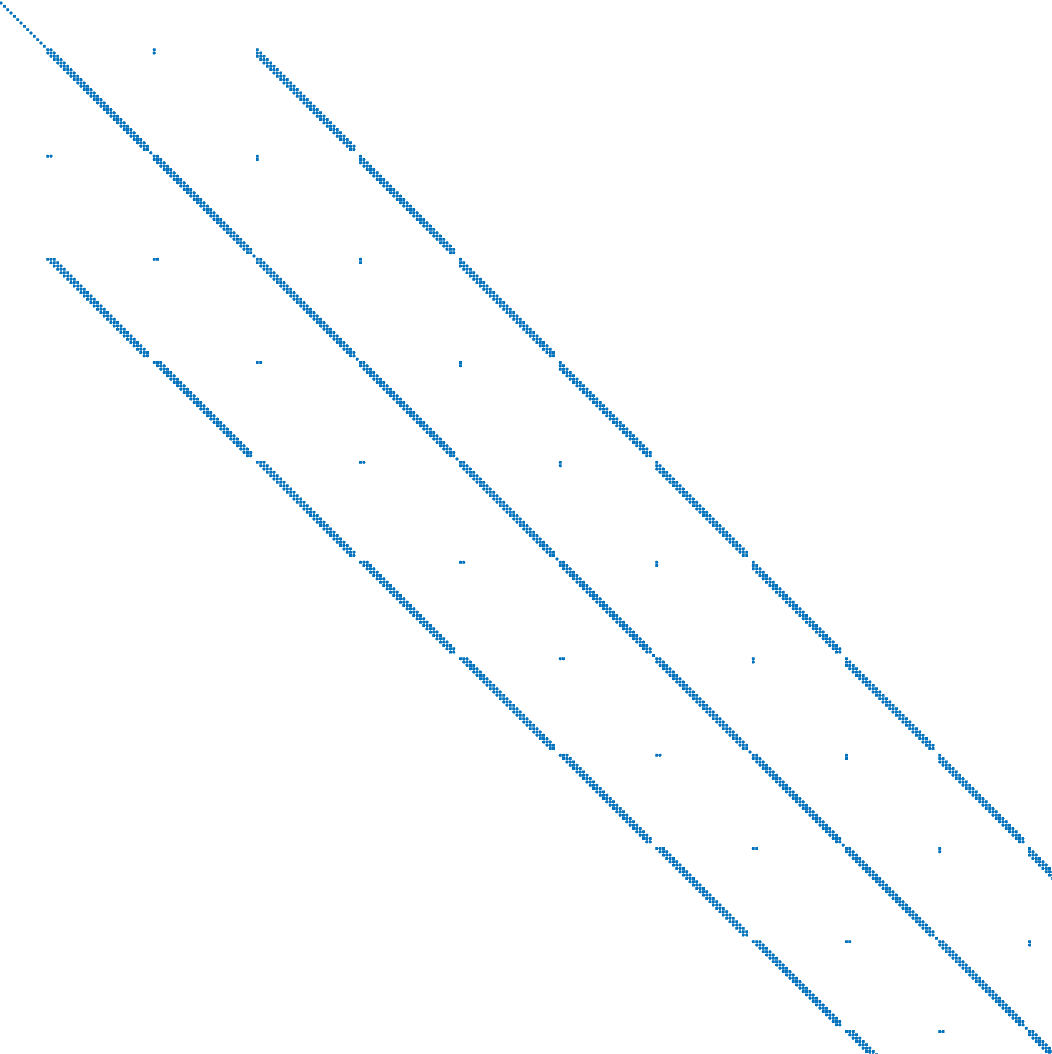}
        \label{fig:quad-sparsity}
    }
    \subfloat[Voronoi\\Sparsity]{
        \includegraphics[width=0.23\textwidth]{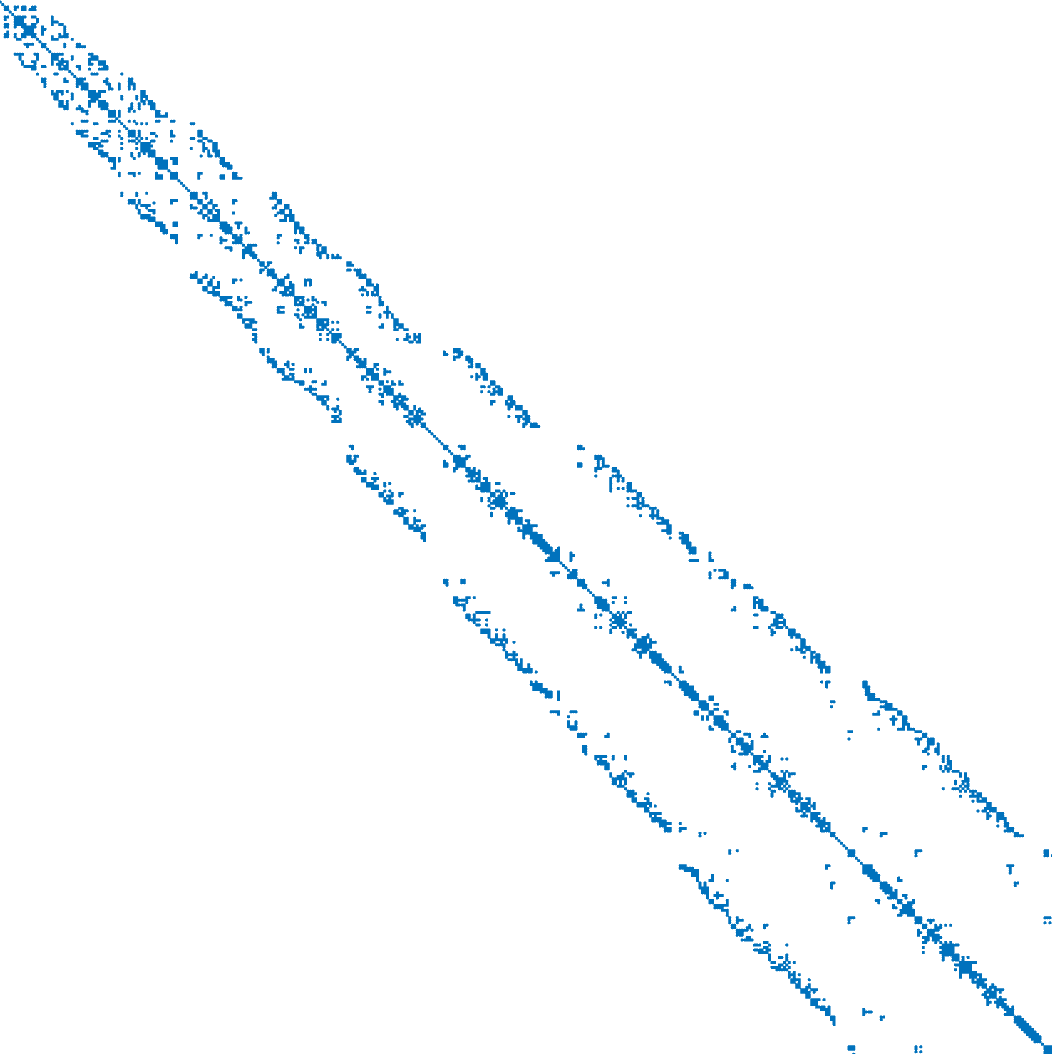}
        \label{fig:voro-sparsity}
    }
    \subfloat[Hexagonal\\Sparsity]{
        \includegraphics[width=0.23\textwidth]{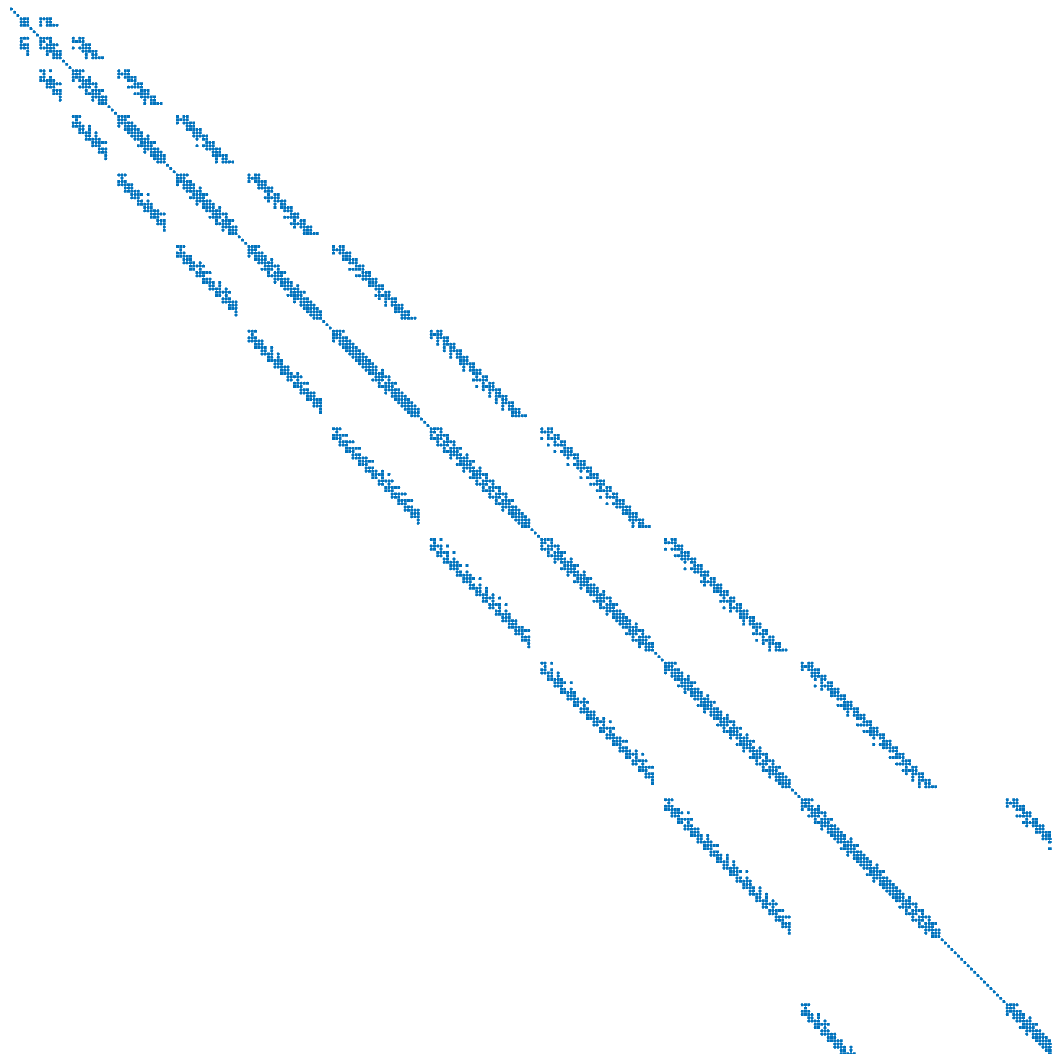}
        \label{fig:hex-sparsity}
    } 
    \caption{The four mesh types considered for two-dimensional test cases (top row). Corresponding sparsity pattern (bottom row) for the diffusion problem discretized with VEM of polynomial degree $p=1$.}
    \label{fig:mesh-types-summary}
\end{figure}

\begin{figure}[!htbp]
    \centering 
    \includegraphics[width=0.9\textwidth]{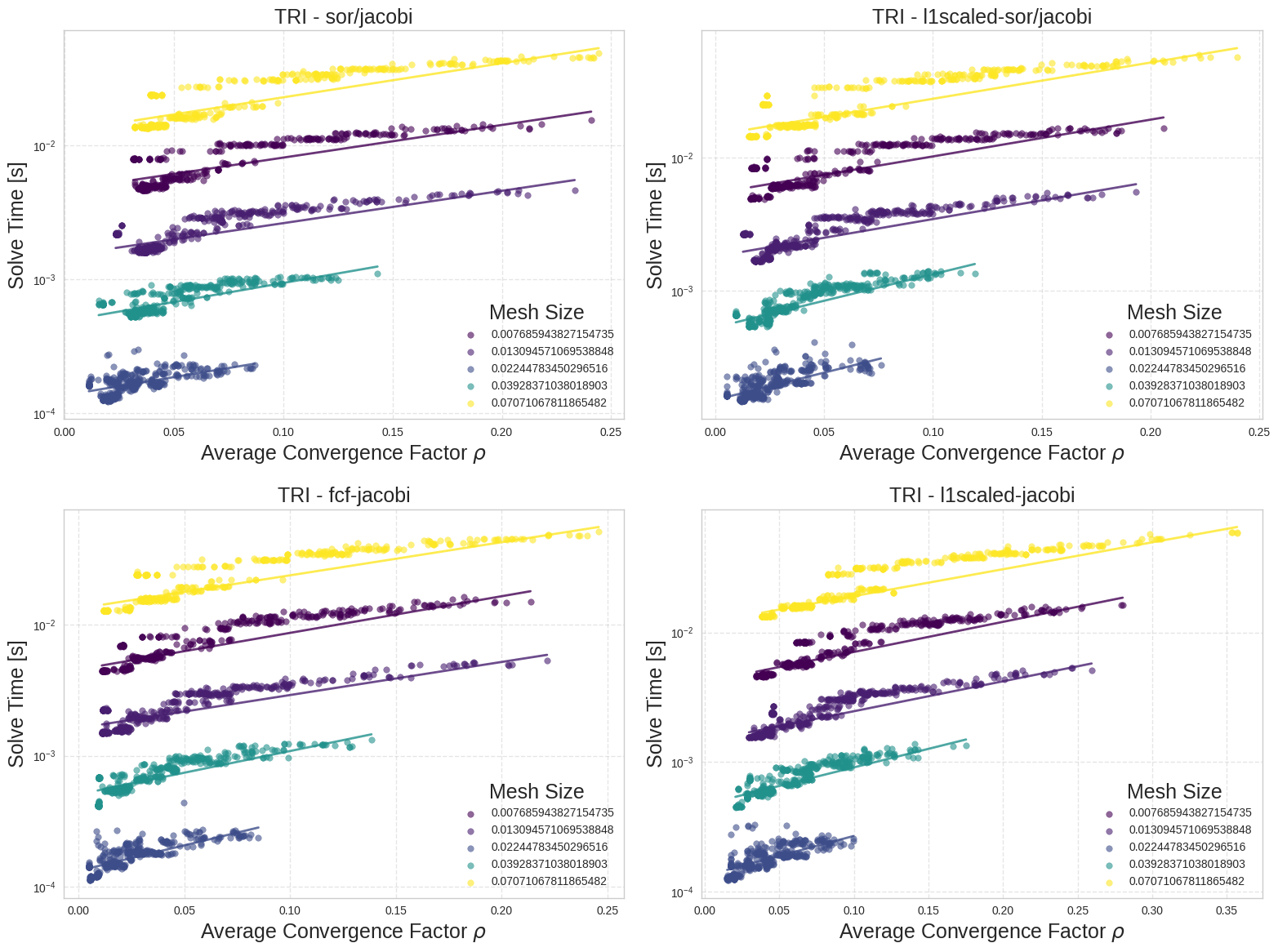}
    \caption{Correlation between the average convergence factor $\rho$ and wall-clock solve time for triangular grids for different smoothers. A clear positive correlation is observed, validating $\rho$ as a reliable proxy for the computational cost. Similar results hold for the other mesh families.}
    \label{fig:rho_vs_time_all}
\end{figure}

This first test case involves a VEM discretization of problem (\ref{eq:model_problem_diffusion}) with discontinuous diffusion coefficients. We remark that this represents a moderately challenging scenario where the default AMG solver is generally effective.
\begin{figure}[!htbp]
    \centering
    \subfloat[Pattern 1 \label{fig:pattern1}]{
\includegraphics[width=0.17\textwidth]{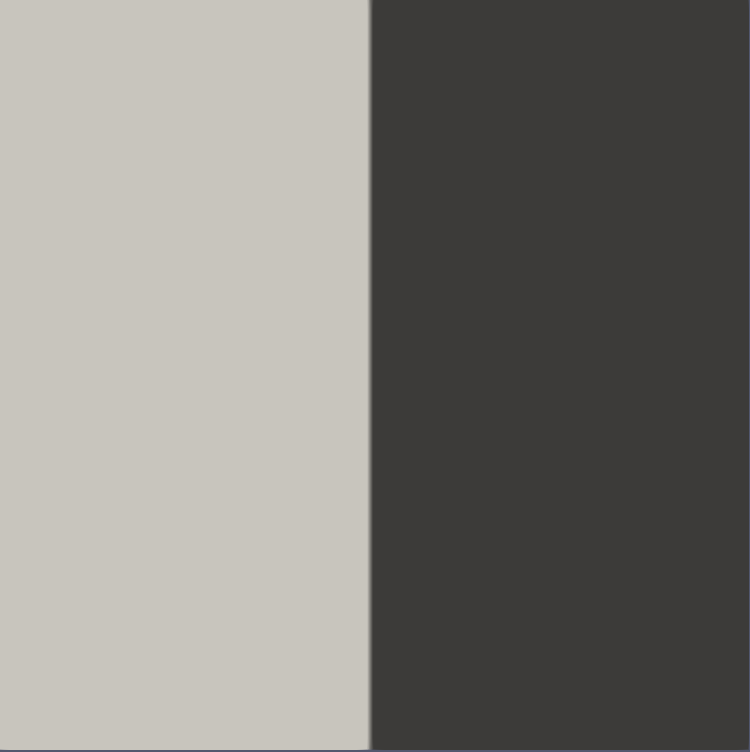}
        }
        \hspace{5mm}
    \subfloat[Pattern 2 \label{fig:pattern2}]{        \includegraphics[width=0.17\textwidth]{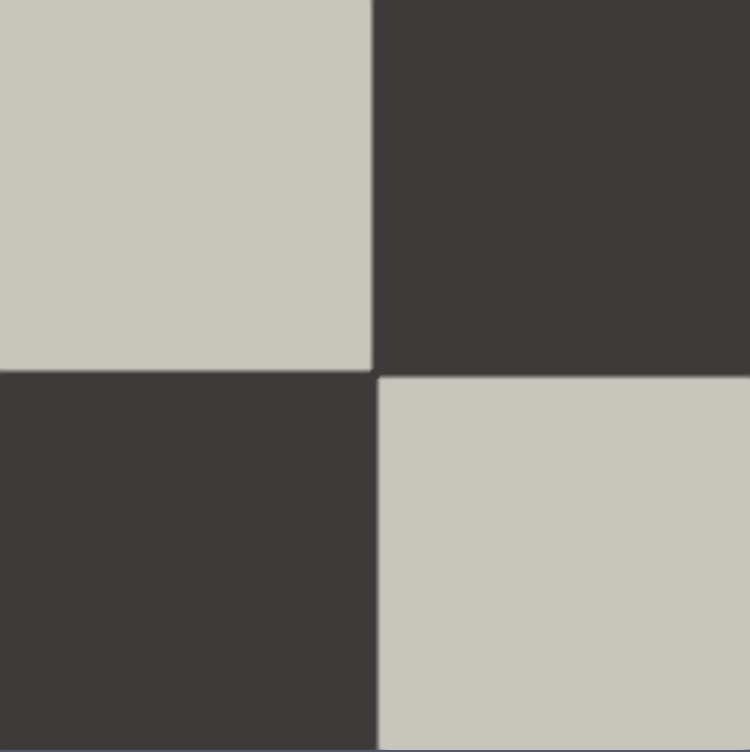}}
        \hspace{5mm}
    \subfloat[Pattern 3 \label{fig:pattern3}]{
\includegraphics[width=0.17\textwidth]{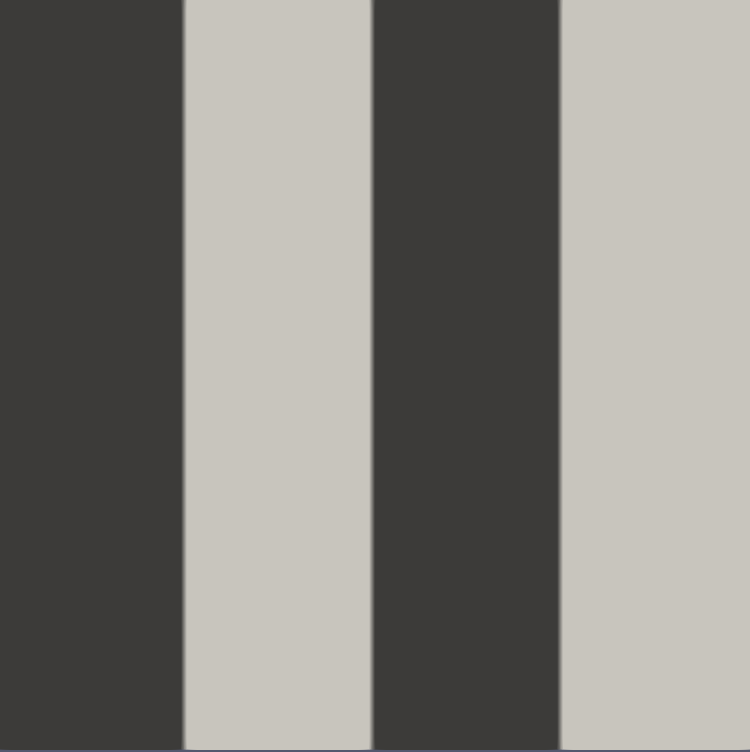}}
        \hspace{5mm}
    \subfloat[Pattern 4 \label{fig:pattern4}]{
\includegraphics[width=0.17\textwidth]{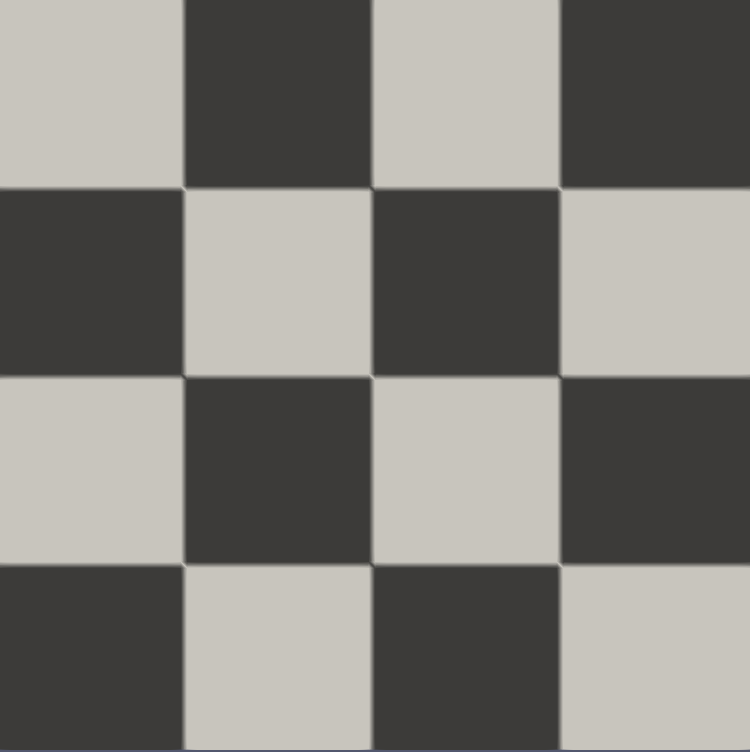}}
    \caption[]{Patterns used to partition the domain $\Omega$ into $\Omega_{white}$ and $\Omega_{gray}$.}
    \label{fig:patterns}
\end{figure}
To define the diffusion coefficient $\kappa$, we partition the domain $\Omega$ into two subdomains $\Omega_{white}$ and $\Omega_{gray}$ according to 4 different patterns (see Figure~\ref{fig:patterns}). In each subdomain we assign a constant value to $\kappa$ in the following way:
\begin{equation}
    \kappa(x,y) = 
    \begin{cases}
        1 & \text{if } (x,y) \in \Omega_{white}, \\
        10^{\epsilon} & \text{if } (x,y) \in \Omega_{gray}.
    \end{cases}
\end{equation}
This diffusion coefficient is discontinuous across subdomains, and the magnitude of the jump is controlled by the parameter $\epsilon \in \mathbb{R}$, which will be varied in the dataset generation.
The domain $\Omega$ is discretized as shown in Figure~\ref{fig:mesh-types-summary}. Namely, we employ 24 distinct meshes, comprising six different refinement levels for each of the four mesh types. The mesh was refined such that the number of elements increases geometrically with each refinement level.\\
In Figure~\ref{fig:rho_vs_time_all}, we plot the correlation between the average convergence factor $\rho$ and the wall-clock solve time for the sequence  of triangular grids and for different smoothers (SOR-Jacobi, $\ell^1$-Jacobi, $\ell^1$-SOR Jacobi and FCF-Jacobi). We observe that $\rho$ and wall clock time are linearly correlated. Similar results were obtained for the other mesh families. Thus, we employ $\rho$ as a measure of the computational cost. 

\begin{figure}[!htbp]
    \centering
    \includegraphics[width=0.45\textwidth]{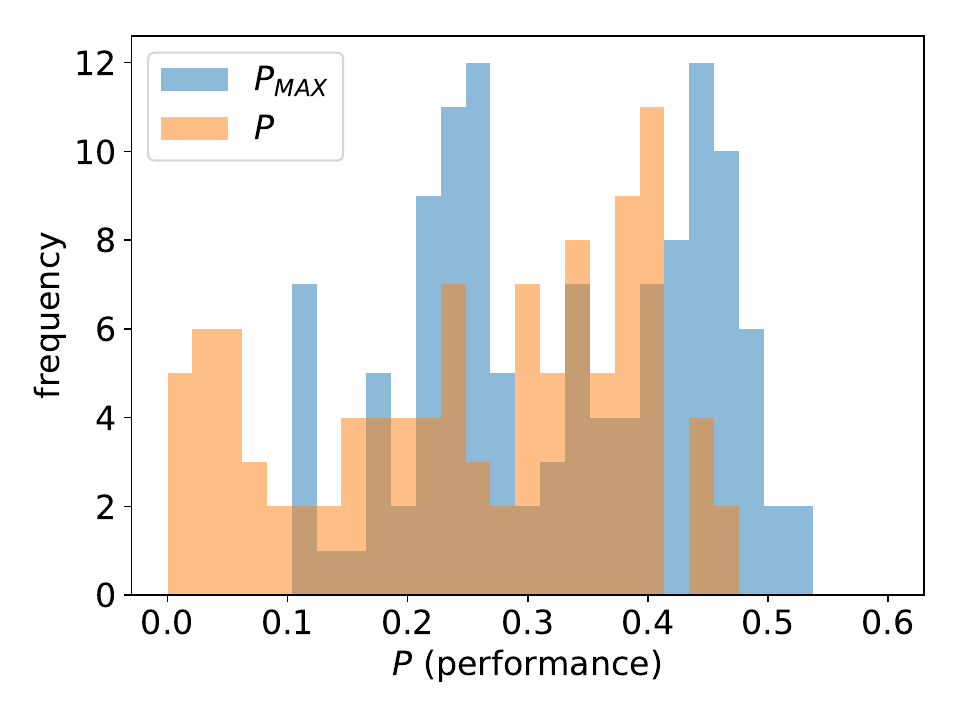}
    \includegraphics[width=0.45\textwidth]{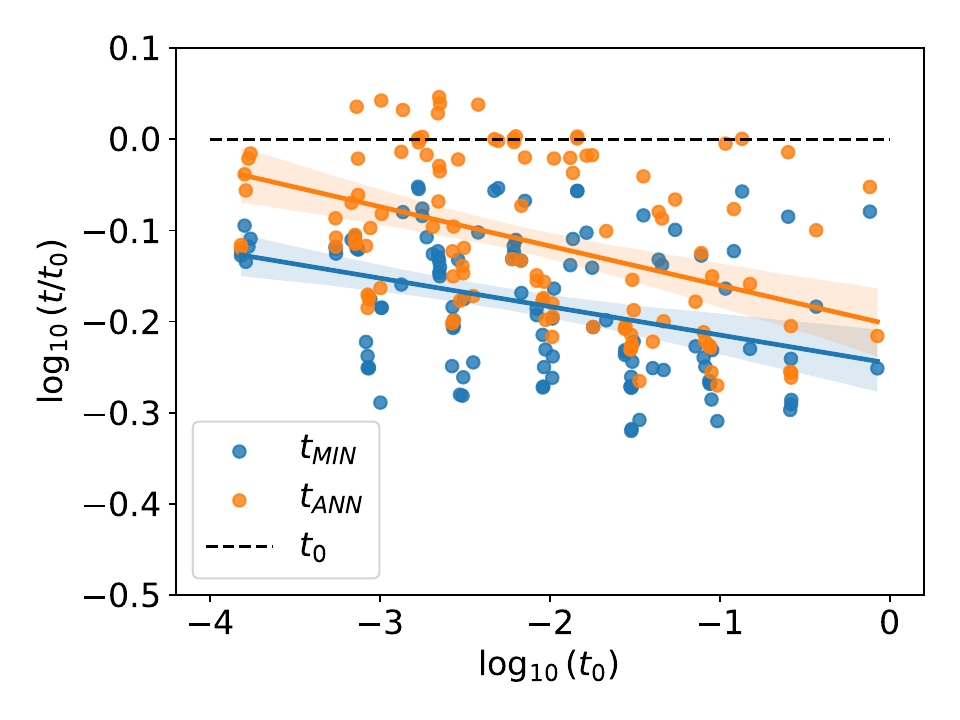}
    \caption{Test Case 1: \textit{Left}: Performance gain $P$ of the AMG-ANN algorithm and maximum theoretical performance $P_{MAX}$. \textit{Right}: Scaling of the AMG-ANN computational cost $t_{ANN}$ and the theoretical minimum $t_{MIN}$ with respect to the default one $t_0$. }
    \label{fig:tcase1}
\end{figure}

In Figure~\ref{fig:tcase1}, we show the gain in performance and the scaling obtained by using our AMG-ANN algorithm. The histogram shows that the gains are evenly distributed, with no peaks, proving that the algorithm is consistent across different scenarios.  The scaling of the computational cost further supports the fact that our algorithm works well as the size of the problem increases, which is of utmost importance for large-scale solvers.\\

The results in Table~\ref{tab:performace} show that the default solver (SOR-Jacobi)  always converges ($P_w=0\%$). Nevertheless, our AMG-ANN algorithm yields substantial improvements. 
The median performance gain $P_M=23.7\%$ indicates a significant reduction in computational time. The algorithm is effective in the vast majority of instances, as shown by $P_B = 87.5\%$.  The performance ratio $P_r=70.1\%$ reveals that the ANN successfully captures over two-thirds of the maximum achievable performance improvement identified in our dataset. This confirms the efficacy of the network's predictions. Finally, the gain in wall-clock time is achieved even though the network was trained using the approximate convergence factor $\rho$. This indicates that our choice of $\rho$ seems an effective and cheaper proxy for computational cost.
\begin{table}[!htbp]
\centering
\begin{tabular}{l|lllllllll}
            & $P_B$ & $P_w$ & $P_m$ & $P_M$ & $P_{MAX}$ & $P_r$ & $P_M^2$ & $P_M^3$ & $P_M^4$ \\
            \hline
Test Case 1 & 87.5  & 0.0     & 21.1  & 23.7  & 33.8      & 70.1  & 17.8    & 22.6    & 19.1    \\
Test Case 2 & 100   & 100   & 100   & 100   & 100       & 100   & 21.9    & 99.6    & 99.7    \\
Test Case 3 & 100   & 100   & 100   & 100   & 100       & 100   & 26.9    & 97.1    & 96.5    \\
Test Case 4 & 100   & 100   & 100   & 100   & 100       & 100   & 24.2    & 100     & 100     \\
Test Case 5 & 100   & 100   & 100   & 100   & 100       & 100   & 24.8    & 97.2    & 95.5    
\end{tabular}
\caption{Evaluation of the performance of the AMG-ANN algorithm for each test case.}
\label{tab:performace}
\end{table}

\subsubsection{Test Case 2: PolyDG Discretization of the Diffusion Problem in 2D}
\begin{figure}[!htbp]
    \centering
    \subfloat[Test Case 2]{
        \includegraphics[width=0.49\textwidth]{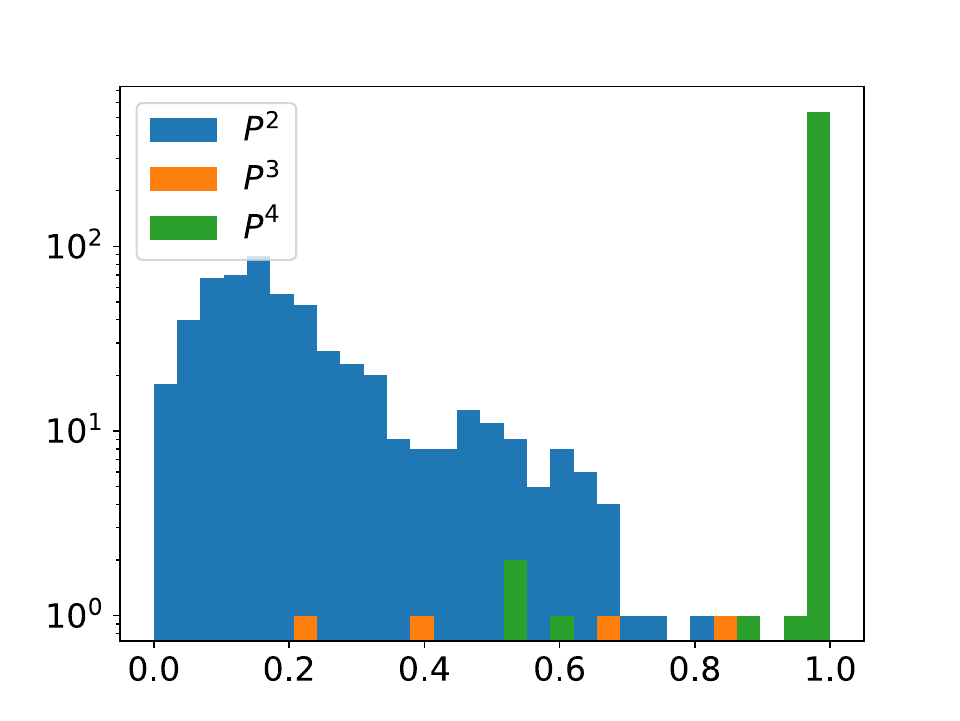}
    }
    \subfloat[Test Case 3]{
        \includegraphics[width=0.49\textwidth]{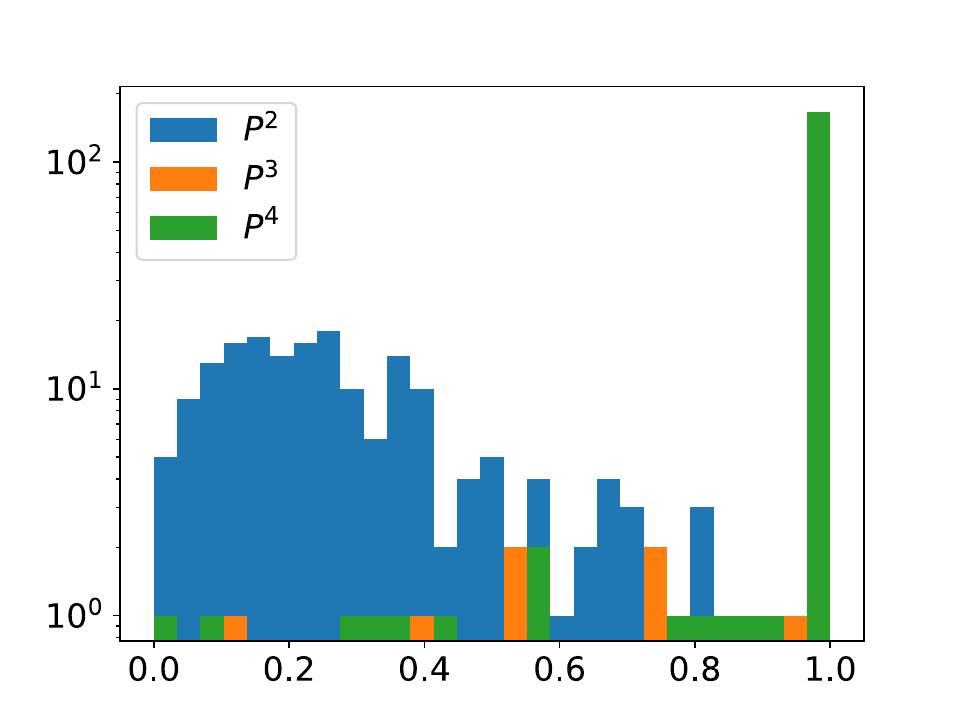}
    }
    \\
    \subfloat[Test Case 4]{
        \includegraphics[width=0.49\textwidth]{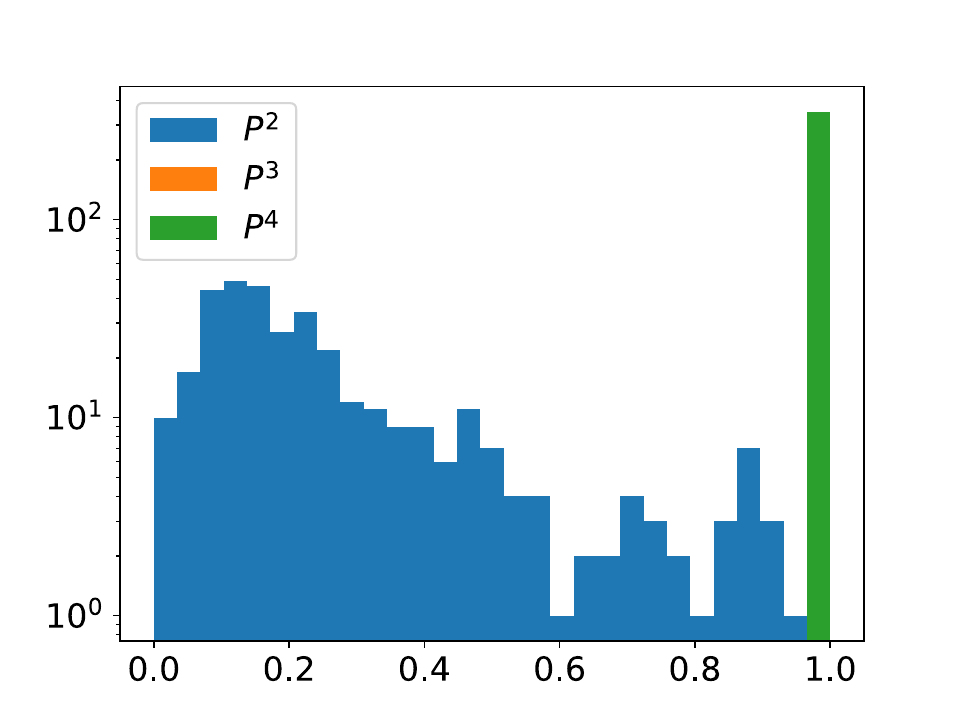}
    }
    \subfloat[Test Case 5]{
     \includegraphics[width=0.49\textwidth]{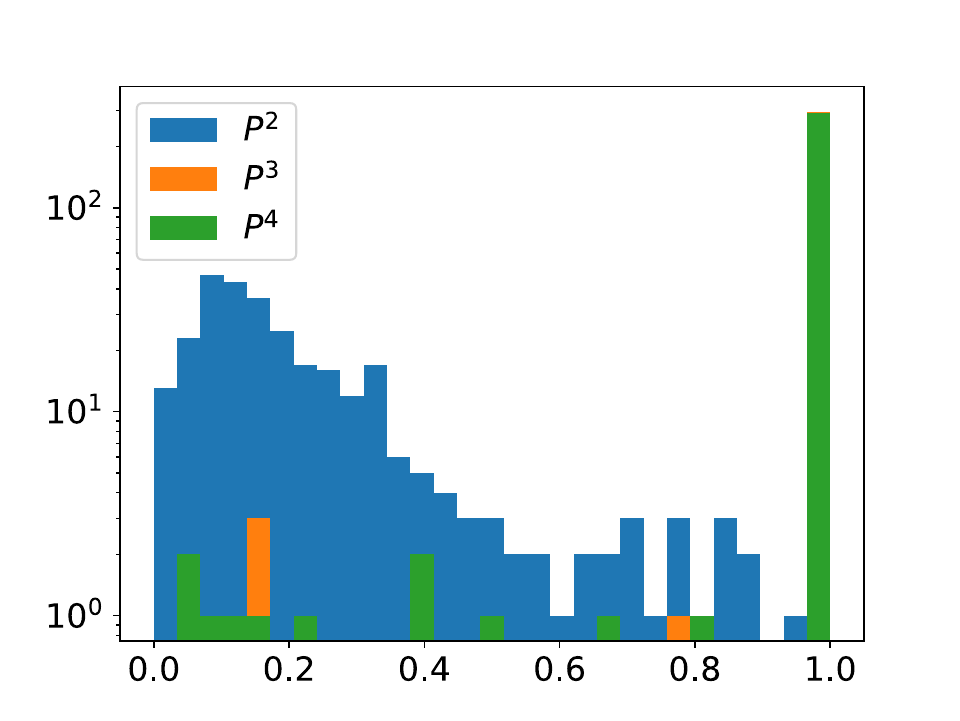}
    }
    \caption{Test Cases 2-5:  Performance gain $P^i$ of the AMG-ANN algorithm with respect to the choice of the smoother: (2) $\ell^1$-Jacobi, (3) $\ell^1$-SOR Jacobi and (4) FCF-Jacobi. The $\ell^1$-SOR Jacobi is usually the best, while the other most often do not reach convergence. Results for (1) SOR-Jacobi are omitted because $P^1=1$ for all test samples. Indeed, when SOR-Jacobi is employed as a smoother, the AMG method fails to converge (see Table~\ref{tab:performace}).} 
    \label{fig:tcase25}
\end{figure}

We consider the diffusion problem \eqref{eq:model_problem_diffusion}, discretized by means of PolyDG with polynomial degree $p=1, 2, 3, 4$. The domain $\Omega=(0, 1)^2$ is a unit square discretized in four different ways, as shown in the top row of Figure~\ref{fig:quad-mesh}. Different levels of refinement are employed when building the dataset. Each refinement level contains a number of elements that grows geometrically relative to the previous level up to $n=200,000$. The diffusion coefficient $\kappa$ is a strongly heterogeneous, piecewise constant and conforming to the mesh. Namely, $\kappa$ takes the form of $\kappa=10^{\epsilon_i}$ on the $i$-th cell of the mesh. The values of $\epsilon_i$ are chosen by extracting values for a uniform distribution in $[0, \epsilon_{MAX}]$, where $\epsilon_{MAX} \in \{1, 2, 4\}$. The penalty $\gamma$ is chosen uniformly at random in $[5, 20]$, making sure that the chosen value guarantees that the corresponding stiffness matrix $A$ is positive definite.
This test case significantly increases the difficulty by employing a PolyDG discretization with highly heterogeneous coefficients. Indeed, the presence of ``duplicate'' degrees of freedom (typical of discontinuous discretizations)  makes the application of the AMG method more challenging. Moreover,  the presence of highly heterogeneous coefficients  produces more ill-conditioned matrices.
The details about the performance are shown in Figure~\ref{fig:tcase25} and Table~\ref{tab:performace}. More precisely, here we address the performance gain $P^i$, $i=2,3,4$, of the AMG-ANN algorithm with respect to the choice of the smoother. Results for smoother (1) (SOR-Jacobi) are omitted because $P^1=1$ for all test samples. Indeed, when SOR-Jacobi is employed as a smoother, the AMG method fails to converge (see Table~\ref{tab:performace}). The smoother (2) is the $\ell^1$-Jacobi methods, whereas the smoothers (3) and (4) are the $\ell^1$-SOR Jacobi and the FCF-Jacobi relaxations, respectively.
The results underscore a critical strength of our approach: its ability to render intractable problems solvable. Indeed, as shown in Table~\ref{tab:performace}, the default AMG configuration fails to converge in every instance ($P_w=100\%$). Consequently, any convergent parameter set represents an (infinite) improvement, leading to perfect scores in $P_B, P_m, P_M, P_{MAX},$ and $P_r$.
The primary contribution of the ANN here is identifying a suitable smoother, as the default SOR-Jacobi method is inadequate. The metrics $P_M^i$ allows to further analyze the algorithm's performance. The near-perfect scores for $P_M^3=99.6\%$ and $P_M^4=99.7\%$ confirm that smoothers (3) $\ell_1$-SOR Jacobi and (4) FCF-Jacobi are poor choices. The much lower value of $P_M^2=21.9\%$ allows to draw two conclusions. First, the smoother (2) $\ell^1$-Jacobi seems to be the most robust choice for this class of problems.
Second, even after making this correct choice, the ANN’s fine-tuning of the threshold $\theta$ provides an additional median time reduction of $21.9\%$. This highlights the dual benefit of our approach: robust smoother selection and effective parameter optimization.

\subsubsection{Test Case 3: PolyDG Discretization of the Diffusion Problem in 3D}
\begin{figure}[!htbp]
    \centering
    \subfloat[Cartesian ]{
        \includegraphics[width=0.2\textwidth]{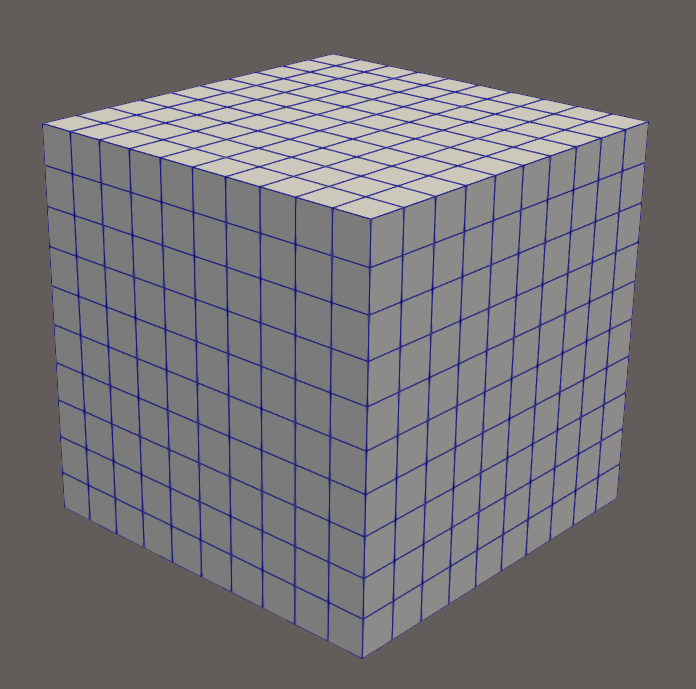}
        }
        \hspace{5mm}
    \subfloat[Structured\\ Tetrahedra ]{
        \includegraphics[width=0.2\textwidth]{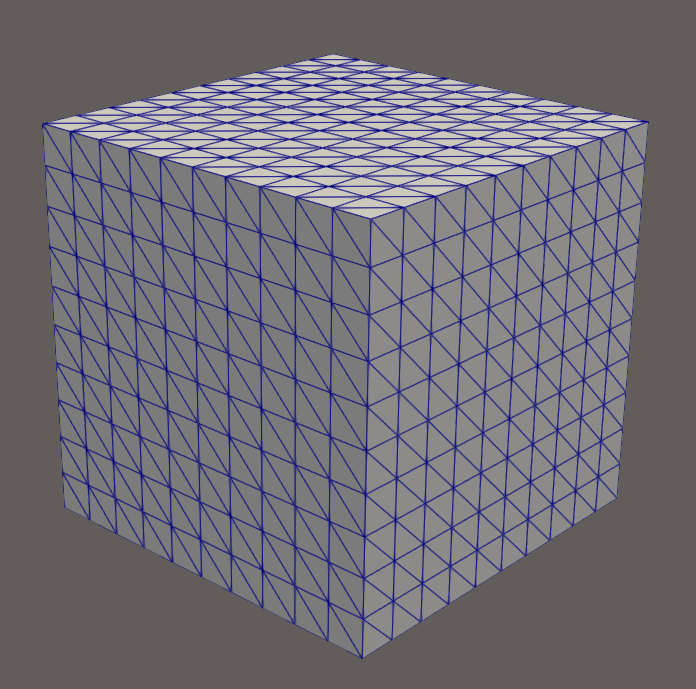}}
        \hspace{5mm}
    \subfloat[Extruded\\ Voronoi ]{
        \includegraphics[width=0.2\textwidth]{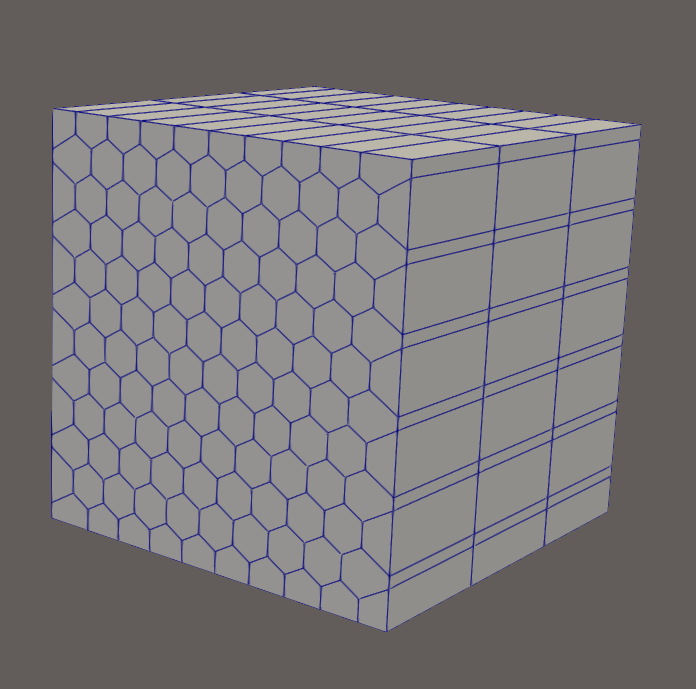}}
        \hspace{5mm}
    \subfloat[Agglomerated]{
        \includegraphics[width=0.2\textwidth]{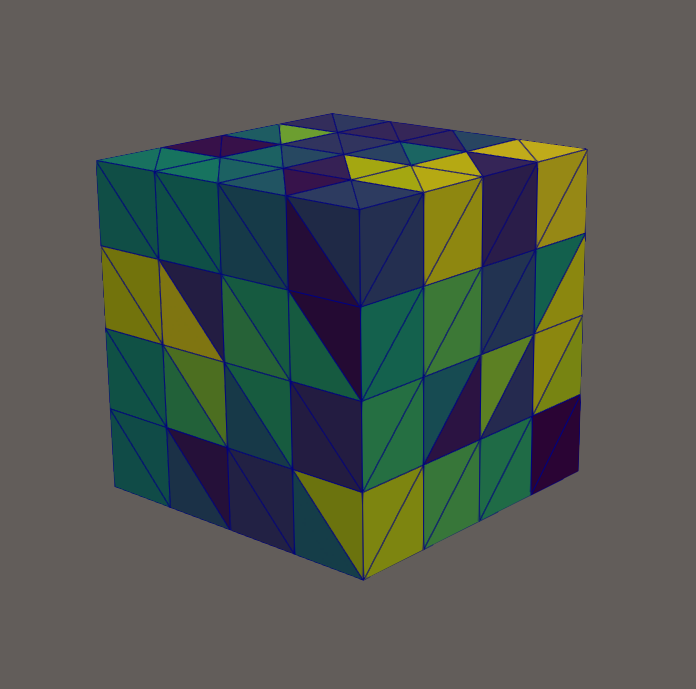}}
    \caption[]{The four mesh types considered in 3D test cases.}
    \label{fig:grids3d}
\end{figure}

We extend the previous test case by considering the three-dimensional version. Namely, we consider  $\Omega = (0, 1)^3$, discretized with four different types of grids, as shown in Figure~\ref{fig:grids3d}. 
We consider various refinements of the grids until we reach $n=100,000$ degrees of freedom.  The diffusion coefficient $\kappa$ takes the same form as the previous test case. Extending the previous scenario to three dimensions further increases the conditioning challenges.
We show the performance of our algorithm in Figure~\ref{fig:tcase25} and Table~\ref{tab:performace} (Test Case 3).
As in two-dimensions, the default AMG solver consistently fails to converge ($P_w=100\%$), making the AMG-ANN's ability to find a convergent set up essential. The results in Table~\ref{tab:performace} show that the landscape of optimal smoothers is more complex in 3D. The performance gains relative to smoothers (3) and (4) seems to remain very high ($P_M^3=97.1\%$ and $P_M^4=96.5\%$), but they are slightly lower than in the 2D case. This seems to suggest that while smoother (2) is still generally the best, the other smoothers are not as uniformly suboptimal. The best choice may therefore be more dependent on the specific matrix properties. The ANN successfully navigates this more complex decision space. However, the benefit of fine-tuning the threshold parameter $\theta$ for the best smoother family is even more pronounced in 3D. This yields a median performance gain of $P_M^2 = 26.9\%$.

\subsection{Numerical Results: Linear Elasticity Problem}\label{sec:numerical-results-elasticity}
We now shift to the discretization of the linear elasticity problem, which introduces further complexities due to its block-structured nature. 

\subsubsection{Test Case 4: PolyDG Discretization of the Linear Elasticity Problem in 2D}
We consider the PolyDG discretization of the linear elasticity problem \ref{eq:model_problem_diffusion}.
The domain $\Omega = (0, 1)^2$ is a unit square discretized with the grids outlined in the first row of Figure~\ref{fig:quad-mesh}. For each kind of grid, we employ different refinements until the total number of degrees of freedom exceeds $n>100,000$. 
To solve the elasticity problem \eqref{eq:model_problem_elasticity}, we consider formulation \eqref{eq:DG_elasticity} and polynomial approximation degrees $p=1,2,3,4$.
We write the Lam\'e parameters in terms of the Young's modulus $E > 0$
and the diffusion ratio $\nu \in (0, \frac{1}{2})$ (see Eq.~\eqref{eq:lame-young}). We fix $\nu = 0.29$ and choose a strongly heterogeneous Young modulus $E$. Namely, for each cell, we set  $E=10^\epsilon$, where $\epsilon$ is chosen randomly in $[0, \epsilon_{MAX}]$, with $\epsilon_{MAX} \in \{1, 2, 4\}$. To enhance the dataset, the penalty coefficient $\gamma$ is chosen in the interval $[4, 20]$. In all the cases the matrix $A$ is symmetric and positive definite.\\

Similar to the diffusion cases, the default AMG parameters frequently fail to yield a convergent solver, making the AMG-ANN indispensable. We have run the same set of experiments as before and the results are reported in Figure~\ref{fig:tcase25} and Table~\ref{tab:performace} (Test Case 4).  From the computed metrics, it seems that the $\ell^1$-scaled SOR Jacobi smoother consistently emerges as the most effective choice across different mesh types and polynomial orders. Our algorithm reliably identifies the optimal smoother in all instances. The primary benefit is therefore transforming a non-convergent method into a robust and efficient one. The median performance gain attributed solely to the optimization of the threshold $\theta$ after selecting the right smoother family is $P_M^2=24.2\%$. This seems to indicate that a careful parameter tuning remains crucial for improving substantially the solver's performance.

\subsubsection{Test Case 5: PolyDG Discretization of the Linear Elasticity Problem in 3D}
The final test case represents the most complex scenario. It combines the challenges of a three-dimensional setting with the vector-valued structure of linear elasticity \eqref{eq:model_problem_elasticity}, and a discontinuous discretization \eqref{eq:DG_elasticity}. The computational domain is $\Omega = (0,1)^3$. As before, a sample of the mesh families considered is shown in Figure~\ref{fig:grids3d}. The meshes are progressively refined until the total number of degrees of freedom $n>10^5$. 
We fix $\nu = 0.29$, while the Young's modulus $E$ is chosen to be strongly heterogeneous across elements. For each cell, we set
\begin{equation}
    E = 10^{\epsilon}, \qquad \epsilon \sim \mathcal{U}(0, \epsilon_{MAX}), 
    \quad \epsilon_{MAX} \in \{1, 2, 4\}.
\end{equation}
The stabilization parameter $\gamma$ is randomly drawn from the interval $[5,20]$, ensuring that the global stiffness matrix $A$ remains symmetric positive definite. The discretization is performed with a PolyDG method of order $p=1,2,3,4$, cf. \eqref{eq:DGspace}. This yields a block-structured stiffness matrices with significant heterogeneity and conditioning challenges.
The numerical results for this case are summarized in Figure~\ref{fig:tcase25} and Table~\ref{tab:performace} (Test Case 5).
From the obtained results,  $\ell^1$-scaled SOR Jacobi is consistently selected by the AMG-ANN algorithm as the most effective smoother. Moreover, it seems robust independently from the discretization order $p$ and the type of mesh.
The baseline AMG setup frequently fails to converge within the prescribed tolerance, in particular for high-order discretizations ($p \geq 3$) and on agglomerated meshes. Our algorithm seems to be able to select in all cases a convergent configuration.
From the numerical results it seems that performance improvements remain substantial. Indeed,  optimizing only the threshold parameter $\theta$ yields a reduction of $24.8\%$ in computational cost.
As opposite to what observed in the corresponding two-dimensional case, the interplay between the smoother choice and threshold parameter parameter tuning is more delicate in 3D. Indeed, some of the smoothers appear to be optimal only for specific mesh families. This can be observed in the degradation of performance in $P_M^3$ and $P_M^4$ relative to $P_M^2$.This can be explain taking into account the higher complexity of the 3D setting.
Overall, these results seem to demonstrate that the proposed AMG-ANN framework scales robustly to heterogeneous 3D linear elasticity discretized by high-order PolyDG methods.

\section{Conclusions}\label{sec:conclusions}
In this work, we employed deep learning to accelerate Algebraic Multigrid iterative solvers for linear systems arising from polytopal discretizations of PDEs.
To improve the performance of AMG solvers, which are highly sensitive to the choice of smoother and parameters, we introduce a novel ANN-AMG framework. ANN-AMG automatically tunes both the smoother selection and the strong threshold parameter in AMG solvers using a novel neural network–driven approach.
By interpreting the matrix as an image and designing a tailored convolutional architecture and a novel pooling strategy, our ANN-AMG achieves robust, on-the-fly parameter selection, eliminating the need for manual fine-tuning. 
The proposed approach is entirely non-intrusive and compatible with existing PDE solvers and parallel AMG implementations. Numerical experiments carried out on both PolyDG and VEM discretizations of diffusion and elasticity problems in two- and three-dimensions demonstrated that we can reduce the computational cost by up to $27\%$ compared to AMG approaches with ``default" choices.
Our results demonstrate the advantages of enhancing algebraic iterative solvers with deep learning to achieve efficient numerical solutions of differential problems posed on complex and heterogeneous domains via polytopal discretizations. 

Further research aims to integrate adaptivity into multilevel solvers, including data-driven coarsening and design of intergrid operators. In addition, we plan to extend the ANN-AMG framework to non-symmetric or indefinite systems, and to explore transferability across broader classes of PDE models and discretizations. Possible future research could also extend the proposed framework to nonlinear algebraic multigrid methods to accelerate convergence for large-scale nonlinear PDE systems, including nonlinear elasticity, multiphase flow, porous media, and reaction–diffusion problems.

\section*{Acknowledgments}
\thanks{This work received funding from the European Union (ERC SyG, NEMESIS, project number 101115663).
Views and opinions expressed are, however, those of the authors only and do not necessarily reflect
those of the European Union or the European Research Council Executive Agency. Neither the
Neither the European Union nor the granting authority can be held responsible for them. PFA, MC, and MV are members of INdAM-GNCS. The present research is part of the activities of “Dipartimento di Eccellenza 2023-2027”, funded by MUR, Italy.
PA and MV also acknowledge MUR--PRIN/PNRR 2022 grant n. P2022BH5CB, funded by MUR.}

\printbibliography

@article{BELLOMO2024,
title = {Life and self-organization on the way to artificial intelligence for collective dynamics},
journal = {Physics of Life Reviews},
volume = {51},
pages = {1-8},
year = {2024},
issn = {1571-0645},
author = {Nicola Bellomo and Marina Dolfin and Jie Liao},
}

@ARTICLE{Quarteroni2025,
	author = {Quarteroni, Alfio and Gervasio, Paola and Regazzoni, Francesco},
	title = {Combining physics-based and data-driven models: advancing the frontiers of research with scientific machine learning},
	year = {2025},
	journal = {Mathematical Models and Methods in Applied Sciences},
	volume = {35},
	number = {4},
	pages = {905 – 1071},
}

@article{BBM_acta_2023,
title={The {V}irtual {E}lement method}, volume={32}, 
journal={Acta Numerica}, 
author={Beir{\~a}o da Veiga, Lourenço and Brezzi, Franco and Marini, L. Donatella and Russo, Alessandro}, year={2023}, 
pages={123–202},
}

@article{AHMAD2013,
title = {Equivalent projectors for {V}irtual {E}lement methods},
journal = {Computers $\&$ Mathematics with Applications},
volume = {66},
number = {3},
pages = {376-391},
year = {2013},
issn = {0898-1221},
author = {B. Ahmad and A. Alsaedi and F. Brezzi and L.D. Marini and A. Russo},
}

@article{BBM_2013,
author = {Beir{\~a}o da Veiga, L.  and Brezzi, F. and Marini, L. D.},
title = {{V}irtual {E}lements for Linear Elasticity Problems},
journal = {SIAM Journal on Numerical Analysis},
volume = {51},
number = {2},
pages = {794-812},
year = {2013},
}

@book{BruntonKutz2019, 
place={Cambridge}, 
title={Data-Driven Science and Engineering: Machine Learning, Dynamical Systems, and Control}, 
publisher={Cambridge University Press}, 
author={Brunton, Steven L. and Kutz, J. Nathan}, year={2019},
}

@ARTICLE{BruntonKutz2018,
	author = {Lusch, Bethany and Kutz, Jose Nathan and Brunton, Steven L.},
	title = {Deep learning for universal linear embeddings of nonlinear dynamics},
	year = {2018},
	journal = {Nature Communications},
	volume = {9},
	number = {1},
}

@ARTICLE{BruntonKutz2016,
	author = {Brunton, Steven L. and Proctor, Joshua L. and Kutz, Jose Nathan},
	title = {Discovering governing equations from data by sparse identification of nonlinear dynamical systems},
	year = {2016},
	journal = {Proceedings of the National Academy of Sciences of the United States of America},
	volume = {113},
	number = {15},
	pages = {3932 - 3937},
}

@ARTICLE{Quarteroni2019,
	author = {Regazzoni, Francesco and Dede', Luca and Quarteroni, Alfio Maria},
	title = {Machine learning for fast and reliable solution of time-dependent differential equations},
	year = {2019},
	journal = {Journal of Computational Physics},
	volume = {397},
	pages = {},
}

@ARTICLE{Zuazua_2023,
	author = {Ruiz-Balet, Domènec and Zuazua, Enrique},
	title = {Neural ODE Control for Classification, Approximation, and Transport},
	year = {2023},
	journal = {SIAM Review},
	volume = {65},
	number = {3},
	pages = {735 - 773},
}

@ARTICLE{Zuazua_2022,
	author = {Geshkovski, Borjan and Zuazua, Enrique},
	title = {{Turnpike in optimal control of PDEs, ResNets, and beyond}},
	year = {2022},
	journal = {Acta Numerica},
	volume = {31},
	pages = {135 - 263},
}

@article{Raissi2019,
  title={Physics-informed neural networks: A deep learning framework for solving forward and inverse problems involving nonlinear partial differential equations},
  author={Raissi, Maziar and Perdikaris, Paris and Karniadakis, George Em},
  journal={Journal of Computational Physics},
  volume={378},
  pages={686--707},
  year={2019}
}

@article{Karniadakis2021,
  title={Physics-informed machine learning},
  author={Karniadakis, George Em and Kevrekidis, Ioannis G and Lu, Lu and Perdikaris, Paris and Wang, Sifan and Yang, Liu},
  journal={Nature Reviews Physics},
  volume={3},
  number={6},
  pages={422--440},
  year={2021}
}

@ARTICLE{Koumoutsakos2020,
	author = {Brunton, Steven L. and Noack, Bernd R. and Koumoutsakos, Petros},
	title = {Machine Learning for Fluid Mechanics},
	year = {2020},
	journal = {Annual Review of Fluid Mechanics},
	volume = {52},
	pages = {477 – 508},
}

@ARTICLE{Bertozzi2019,
	author = {Bertozzi, Andrea L. and Merkurjev, Ekaterina},
	title = {Graph-based optimization approaches for machine learning, uncertainty quantification and networks},
	year = {2019},
	journal = {Handbook of Numerical Analysis},
	volume = {20},
	pages = {503 – 531},
	url = {https://www.scopus.com/inward/record.uri?eid=2-s2.0-85065410921&doi=10.1016%2fbs.hna.2019.04.001&partnerID=40&md5=b03e9d9f4a364dc8cad15aa7be96747c}
}

@ARTICLE{Bertozzi2012,
	author = {Bertozzi, Andrea L. and Flenner, Arjuna},
	title = {Diffuse interface models on graphs for classification of high dimensional data},
	year = {2012},
	journal = {Multiscale Modeling and Simulation},
	volume = {10},
	number = {3},
	pages = {1090 – 1118},
	url = {https://www.scopus.com/inward/record.uri?eid=2-s2.0-84867006275&doi=10.1137%2f11083109X&partnerID=40&md5=77ff5c8aada1aa3b3f6a816142c4e396}
}

@book{Antonietti2022VEM,
  editor    = {Antonietti, Paola F. and Beir{\~a}o da Veiga, Louren{\c{c}}o and Manzini, Gianmarco},
  title     = {{The {V}irtual {E}lement Method and its Applications}},
  series    = {{SEMA SIMAI Springer Series}},
  volume    = {31},
  publisher = {Springer},
  address   = {Cham, Switzerland},
  year      = {2022},
  isbn      = {978-3-030-95318-8},
  doi       = {10.1007/978-3-030-95319-5},
  url       = {https://link.springer.com/book/10.1007/978-3-030-95319-5}
}

@article{BeiraoDaVeiga_Brezzi_Marini_Russo_2014,
author = {Beir\~{a}o da Veiga, L. and Brezzi, F. and Marini, L. D. and Russo, A.},
title = {The Hitchhiker's Guide to the {V}irtual {E}lement Method},
journal = {Mathematical Models and Methods in Applied Sciences},
volume = {24},
number = {08},
pages = {1541-1573},
year = {2014},
}

@ARTICLE{BeiraoDaVeiga_Lovadina_Mora_2015,
	author = {Beir{\~a}o da Veiga, L. and Lovadina, C. and Mora, D.},
	title = {A {V}irtual {E}lement Method for elastic and inelastic problems on polytope meshes},
	year = {2015},
	journal = {Computer Methods in Applied Mechanics and Engineering},
	volume = {295},
	pages = {327 -- 346},
}

@ARTICLE{BeiraoDaVeiga_Dassi_Russo_2017,
	author = {Beir{\~a}o da Veiga, L. and Dassi, F. and Russo, A.},
	title = {High-order {V}irtual {E}lement Method on polyhedral meshes},
	year = {2017},
	journal = {Computers and Mathematics with Applications},
	volume = {74},
	number = {5},
	pages = {1110 -- 1122},
}

@ARTICLE{BeiraoDaVeiga_Brezzi_Marini_Russo_2016,
	author = {Beir{\~a}o Da Veiga, L. and Brezzi, F. and Marini, L.D. and Russo, A.},
	title = {{V}irtual {E}lement Method for general second-order elliptic problems on polygonal meshes},
	year = {2016},
	journal = {Mathematical Models and Methods in Applied Sciences},
	volume = {26},
	number = {4},
	pages = {729 -- 750},
}

@ARTICLE{BeiraoDaVeiga_Brezzi_Marini_2013,
	author = {Beir{\~a}o Da Veiga, L. and Brezzi, F. and Marini, L.D.},
	title = {{V}irtual {E}lements for linear elasticity problems},
	year = {2013},
	journal = {SIAM Journal on Numerical Analysis},
	volume = {51},
	number = {2},
	pages = {794 -- 812},
}

@ARTICLE{Arnold_Brezzi_Cockburn_Marini_2001,
	author = {Arnold, Douglas N. and Brezzi, Franco and Cockburn, Bernardo and Donatella Marini, L.},
	title = {Unified analysis of discontinuous {Galerkin} methods for elliptic problems},
	year = {2001},
	journal = {SIAM Journal on Numerical Analysis},
	volume = {39},
	number = {5},
	pages = {1749 -- 1779},
}

@ARTICLE{Feder2025,
	author = {Feder, Marco and Cangiani, Andrea and Heltai, Luca},
	title = {{R3MG: R-tree based agglomeration of polytopal grids with applications to multilevel methods}},
	year = {2025},
	journal = {Journal of Computational Physics},
	volume = {526},
}

@ARTICLE{BeiraodaVeiga_Lipnikov_Manzini_2014,
	author = {Beir{\~a}o da Veiga, Lourenço and Lipnikov, Konstantin and Manzini, Gianmarco},
	title = {The mimetic finite difference method for elliptic problems},
	year = {2014},
	journal = {Modeling, Simulation and Applications},
	volume = {11},
	pages = {1 -- 389},
}

@ARTICLE{DiPietroErnLemair_2014,
	author = {Di Pietro, Daniele A. and Ern, Alexandre and Lemaire, Simon},
	title = {An arbitrary-order and compact-stencil discretization of diffusion on general meshes based on local reconstruction operators},
	year = {2014},
	journal = {Computational Methods in Applied Mathematics},
	volume = {14},
	number = {4},
	pages = {461 -- 472},
}

@ARTICLE{Antonietti_Brezzi_Marini_2009,
	author = {Antonietti, Paola F. and Brezzi, Franco and Marini, L. Donatella},
	title = {Bubble stabilization of discontinuous {Galerkin} methods},
	year = {2009},
	journal = {Computer Methods in Applied Mechanics and Engineering},
	volume = {198},
	number = {21-26},
	pages = {1651 -- 1659},
}

@ARTICLE{Bassi_et_al_2012,
	author = {Bassi, F. and Botti, L. and Colombo, A. and Di Pietro, D.A. and Tesini, P.},
	title = {On the flexibility of agglomeration-based physical space discontinuous {Galerkin} discretizations},
	year = {2012},
	journal = {Journal of Computational Physics},
	volume = {231},
	number = {1},
	pages = {45 -- 65},
}

@book{CangianiDongGeorgoulisHouston_2017, 
author = {Cangiani, Andrea and Dong, Zhaonan and Georgoulis, Emmanuil H. and Houston, Paul}, 
title = {hp-Version Discontinuous {Galerkin} Methods on Polygonal and Polyhedral Meshes}, year = {2017}, 
isbn = {3319676717}, 
publisher = {Springer Publishing Company, Incorporated}, 
edition = {1st},
}

@incollection{ruge1987algebraic,
  title={{A}lgebraic {M}ultigrid},
  author={Ruge, John W and St{\"u}ben, Klaus},
  booktitle={{M}ultigrid {M}ethods},
  pages={73--130},
  year={1987},
  publisher={SIAM}
}

@article{yang2002boomeramg,
  title={Boomer{AMG}: A parallel {A}lgebraic {M}ultigrid solver and preconditioner},
  author={Yang, Ulrike Meier and others},
  journal={Applied Numerical Mathematics},
  volume={41},
  number={1},
  pages={155--177},
  year={2002},
  publisher={Elsevier}
}

@inproceedings{falgout2002hypre,
  title={hypre: A library of high performance preconditioners},
  author={Falgout, Robert D and Yang, Ulrike Meier},
  booktitle={International Conference on computational science},
  pages={632--641},
  year={2002},
  organization={Springer}
}

@book{wesseling2004introduction,
  title={An Introduction to {M}ultigrid {M}ethods},
  author={Wesseling, P.},
  isbn={9781930217089},
  lccn={91024430},
  series={An Introduction to Multigrid Methods},
  year={2004},
  publisher={R.T. Edwards}
}

@book{bramble2019multigrid,
  title={{M}ultigrid {M}ethods},
  author={Bramble, James H},
  year={2019},
  publisher={Chapman and Hall/CRC}
}

@book{brandt1983algebraic,
  title={{A}lgebraic {M}ultigrid ({AMG}) for Automatic Multigrid Solution with Application to Geodetic Computations},
  author={Brandt, A. and McCormick, S. and Ruge, J.},
  publisher={National Geodetic Survey and Air Force Office of Scientific Research and National Science Foundation},
  year={1983}
}

@article{brandt1984algebraic,
  title={{A}lgebraic {M}ultigrid ({AMG}) for sparse matrix eqations},
  author={Brandt, Achi},
  journal={Sparsity and its Applications},
  pages={257--284},
  year={1984},
  publisher={Cambridge University Press}
}

@article{zikatanov2008two,
  title={Two-sided bounds on the convergence rate of two-level methods},
  author={Zikatanov, Ludmil T},
  journal={Numerical Linear Algebra with Applications},
  volume={15},
  number={5},
  pages={439--454},
  year={2008},
  publisher={Wiley Online Library}
}

@article{falgout2004generalizing,
  title={On generalizing the {A}lgebraic {M}ultigrid framework},
  author={Falgout, Robert D and Vassilevski, Panayot S},
  journal={SIAM Journal on Numerical Analysis},
  volume={42},
  number={4},
  pages={1669--1693},
  year={2004},
  publisher={SIAM}
}

@article{brezina2006adaptive,
  title={Adaptive {A}lgebraic {M}ultigrid},
  author={Brezina, Marian and Falgout, R and MacLachlan, Scott and Manteuffel, T and McCormick, S and Ruge, John},
  journal={SIAM Journal on Scientific Computing},
  volume={27},
  number={4},
  pages={1261--1286},
  year={2006},
  publisher={SIAM}
}

@article{caldana2024deep,
  title={A deep learning algorithm to accelerate algebraic multigrid methods in finite element solvers of 3{D} elliptic {PDE}s},
  author={Caldana, Matteo and Antonietti, Paola F and others},
  journal={Computers \& Mathematics with Applications},
  volume={167},
  pages={217--231},
  year={2024},
  publisher={Elsevier}
}

@article{brandt1986algebraic,
  title={{A}lgebraic {M}ultigrid theory: The symmetric case},
  author={Brandt, Achi},
  journal={Applied {M}athematics and {C}omputation},
  volume={19},
  number={1-4},
  pages={23--56},
  year={1986},
  publisher={Elsevier}
}

@book{trottenberg2000multigrid,
  title={Multigrid},
  author={Trottenberg, Ulrich and Oosterlee, Cornelius W and Schuller, Anton},
  year={2000},
  publisher={Elsevier}
}

@article{van2001convergence,
  title={Convergence of {A}lgebraic {M}ultigrid based on smoothed aggregation},
  author={Van\v{e}k, Petr and Brezina, Marian and Mandel, Jan},
  journal={Numerische Mathematik},
  volume={88},
  pages={559--579},
  year={2001},
  publisher={Springer}
}

@article{vanek1996algebraic,
  title={{A}lgebraic {M}ultigrid by smoothed aggregation for second and fourth order elliptic problems},
  author={Van\v{e}k, Petr and Mandel, Jan and Brezina, Marian},
  journal={Computing},
  volume={56},
  number={3},
  pages={179--196},
  year={1996},
  publisher={Citeseer}
}

@article{xu2017algebraic,
  title={{A}lgebraic {M}ultigrid methods},
  author={Xu, Jinchao and Zikatanov, Ludmil},
  journal={Acta Numerica},
  volume={26},
  pages={591--721},
  year={2017},
  publisher={Cambridge University Press}
}

@article{falgout2005two,
  title={On two-grid convergence estimates},
  author={Falgout, Robert D and Vassilevski, Panayot S and Zikatanov, Ludmil T},
  journal={Numerical linear algebra with applications},
  volume={12},
  number={5-6},
  pages={471--494},
  year={2005},
  publisher={Wiley Online Library}
}

@article{antonietti2017multigrid,
  title={Multigrid algorithms for hp-version interior penalty {D}iscontinuous {G}alerkin methods on polygonal and polyhedral meshes},
  author={Antonietti, Paola F and Houston, Paul and Hu, Xiaozhe and Sarti, Marco and Verani, Marco},
  journal={Calcolo},
  volume={54},
  pages={1169--1198},
  year={2017},
  publisher={Springer}
}

@article{antonietti2018multigrid,
  title={A multigrid algorithm for the p-version of the {V}irtual {E}lement {M}ethod},
  author={Antonietti, Paola F and Mascotto, Lorenzo and Verani, Marco},
  journal={ESAIM: Mathematical Modelling and Numerical Analysis},
  volume={52},
  number={1},
  pages={337--364},
  year={2018},
  publisher={EDP sciences}
}

@article{pan2022agglomeration,
  title={Agglomeration-based geometric multigrid solvers for compact discontinuous {G}alerkin discretizations on unstructured meshes},
  author={Pan, Yulong and Persson, P-O},
  journal={Journal of Computational Physics},
  volume={449},
  pages={110775},
  year={2022},
  publisher={Elsevier}
}

@article{antonietti2023accelerating,
  title={Accelerating {A}lgebraic {M}ultigrid methods via {A}rtificial {N}eural {N}etworks},
  author={Antonietti, Paola F and Caldana, Matteo and Dede', Luca},
  journal={Vietnam Journal of Mathematics},
  pages={1--36},
  year={2023},
  publisher={Springer}
}

@book{saad2003iterative,
  title={Iterative methods for sparse linear systems},
  author={Saad, Yousef},
  year={2003},
  publisher={SIAM}
}

@article{beirao2014hitchhiker,
  title={The hitchhiker's guide to the {V}irtual {E}lement method},
  author={Beir{\~a}o da Veiga, L and Brezzi, Franco and Marini, Luisa Donatella and Russo, Alessandro},
  journal={Mathematical models and methods in applied sciences},
  volume={24},
  number={08},
  pages={1541--1573},
  year={2014},
  publisher={World Scientific}
}

@article{beirao2013basic,
  title={Basic principles of {V}irtual {E}lement methods},
  author={Beir{\~a}o da Veiga, Lourenco and Brezzi, Franco and Cangiani, Andrea and Manzini, Gianmarco and Marini, L Donatella and Russo, Alessandro},
  journal={Mathematical Models and Methods in Applied Sciences},
  volume={23},
  number={01},
  pages={199--214},
  year={2013},
  publisher={World Scientific}
}

@article{cangiani2014hp,
  title={hp-version discontinuous {G}alerkin methods on polygonal and polyhedral meshes},
  author={Cangiani, Andrea and Georgoulis, Emmanuil H and Houston, Paul},
  journal={Mathematical Models and Methods in Applied Sciences},
  volume={24},
  number={10},
  pages={2009--2041},
  year={2014},
  publisher={World Scientific}
}

@article{antonietti2013hp,
  title={hp-version composite discontinuous {G}alerkin methods for elliptic problems on complicated domains},
  author={Antonietti, Paola F and Giani, Stefano and Houston, Paul},
  journal={SIAM Journal on Scientific Computing},
  volume={35},
  number={3},
  pages={A1417--A1439},
  year={2013},
  publisher={SIAM}
}

@incollection{antonietti2016review,
  title={Review of discontinuous {G}alerkin finite element methods for partial differential equations on complicated domains},
  author={Antonietti, Paola F and Cangiani, Andrea and Collis, Joe and Dong, Zhaonan and Georgoulis, Emmanuil H and Giani, Stefano and Houston, Paul},
  booktitle={Building bridges: connections and challenges in modern approaches to numerical partial differential equations},
  pages={281--310},
  year={2016},
  publisher={Springer}
}

@article{di2020hybrid,
  title={The hybrid high-order method for polytopal meshes},
  author={Di Pietro, Daniele Antonio and Droniou, J{\'e}r{\^o}me},
  journal={Number 19 in Modeling, Simulation and Application},
  volume={84},
  year={2020},
  publisher={Springer}
}

@article{lipnikov2014mimetic,
  title={Mimetic finite difference method},
  author={Lipnikov, Konstantin and Manzini, Gianmarco and Shashkov, Mikhail},
  journal={Journal of Computational Physics},
  volume={257},
  pages={1163--1227},
  year={2014},
  publisher={Elsevier}
}

@inproceedings{cockburn2010hybridizable,
  title={The hybridizable discontinuous {G}alerkin methods},
  author={Cockburn, Bernardo},
  booktitle={Proceedings of the International Congress of Mathematicians 2010 (ICM 2010) (In 4 Volumes) Vol. I: Plenary Lectures and Ceremonies Vols. II--IV: Invited Lectures},
  pages={2749--2775},
  year={2010},
  organization={World Scientific}
}

@article{nguyen2010hybridizable,
  title={A hybridizable discontinuous {G}alerkin method for {S}tokes flow},
  author={Nguyen, Ngoc Cuong and Peraire, Jaime and Cockburn, Bernardo},
  journal={Computer Methods in Applied Mechanics and Engineering},
  volume={199},
  number={9-12},
  pages={582--597},
  year={2010},
  publisher={Elsevier}
}

@article{wang2016weak,
  title={A weak {G}alerkin finite element method for the {S}tokes equations},
  author={Wang, Junping and Ye, Xiu},
  journal={Advances in Computational Mathematics},
  volume={42},
  number={1},
  pages={155--174},
  year={2016},
  publisher={Springer}
}

@article{wang2013weak,
  title={A weak {G}alerkin finite element method for second-order elliptic problems},
  author={Wang, Junping and Ye, Xiu},
  journal={Journal of Computational and Applied Mathematics},
  volume={241},
  pages={103--115},
  year={2013},
  publisher={Elsevier}
}

@ARTICLE{antonietti2025magnet,
	author = {Antonietti, Paola F. and Caldana, Matteo and Mazzieri, Ilario and Fraschini, Andrea Re},
	title = {{MAGNET: an open-source library for mesh agglomeration by {G}raph {N}eural {N}etworks}},
	year = {2025},
	journal = {Engineering with Computers},
	volume = {41},
	number = {6},
	pages = {4825 – 4850}
}

@article{antonietti2024agglomeration,
  title={Agglomeration of polygonal grids using {G}raph {N}eural {N}etworks with applications to multigrid solvers},
  author={Antonietti, Paola F and Farenga, Nicola and Manuzzi, Enrico and Martinelli, Gabriele and Saverio, Luca},
  journal={Computers \& Mathematics with Applications},
  volume={154},
  pages={45--57},
  year={2024},
  publisher={Elsevier}
}

@article{savitzky1964smoothing,
  title={Smoothing and differentiation of data by simplified least squares procedures.},
  author={Savitzky, Abraham and Golay, Marcel JE},
  journal={Analytical Chemistry},
  volume={36},
  number={8},
  pages={1627--1639},
  year={1964},
  publisher={ACS Publications}
}

@article{antonietti2023agglomeration,
  title={Agglomeration-based geometric multigrid schemes for the {V}irtual {E}lement {M}ethod},
  author={Antonietti, Paola F and Berrone, Stefano and Busetto, Martina and Verani, Marco},
  journal={SIAM Journal on Numerical Analysis},
  volume={61},
  number={1},
  pages={223--249},
  year={2023},
  publisher={SIAM}
}

@article{baker2011multigrid,
  title={Multigrid smoothers for ultraparallel computing},
  author={Baker, Allison H and Falgout, Robert D and Kolev, Tzanio V and Yang, Ulrike Meier},
  journal={SIAM Journal on Scientific Computing},
  volume={33},
  number={5},
  pages={2864--2887},
  year={2011},
  publisher={SIAM}
}

@article{hessenthaler2020multilevel,
  title={Multilevel convergence analysis of multigrid-reduction-in-time},
  author={Hessenthaler, Andreas and Southworth, Ben S and Nordsletten, David and Rohrle, Oliver and Falgout, Robert D and Schroder, Jacob B},
  journal={SIAM Journal on Scientific Computing},
  volume={42},
  number={2},
  pages={A771--A796},
  year={2020},
  publisher={SIAM}
}

@inproceedings{akiba2019optuna,
  title={Optuna: A next-generation hyperparameter optimization framework},
  author={Akiba, Takuya and Sano, Shotaro and Yanase, Toshihiko and Ohta, Takeru and Koyama, Masanori},
  booktitle={Proceedings of the 25th ACM SIGKDD international conference on knowledge discovery \& data mining},
  pages={2623--2631},
  year={2019}
}

@inproceedings{loshchilov2017decoupled,
title={Decoupled Weight Decay Regularization},
author={Ilya Loshchilov and Frank Hutter},
booktitle={International Conference on Learning Representations},
year={2019},
url={https://openreview.net/forum?id=Bkg6RiCqY7},
}

@article{dassi2025vem++,
  title={Vem++, a {C}++ library to handle and play with the {V}irtual {E}lement Method},
  author={Dassi, Franco},
  journal={Numerical Algorithms},
  pages={1--43},
  year={2025},
  publisher={Springer}
}

@book{balay2019petsc,
  title={PETSc users manual},
  author={Balay, Satish and Abhyankar, Shrirang and Adams, Mark and Brown, Jed and Brune, Peter and Buschelman, Kris and Dalcin, Lisandro and Dener, Alp and Eijkhout, Victor and Gropp, William and others},
  year={2019},
  publisher={Argonne National Laboratory}
}
\end{document}